\documentclass[10pt]{article}

\RequirePackage{amsmath}
\RequirePackage{amsfonts}
\RequirePackage{amssymb}
\RequirePackage{hyperref}
\RequirePackage{multirow}
\RequirePackage{algorithm}
\RequirePackage{algpseudocode}
\RequirePackage{graphicx}
\RequirePackage{xcolor}
\RequirePackage{amsthm}
\RequirePackage{amsmath}
\RequirePackage{amsfonts}
\RequirePackage{amssymb}
\RequirePackage{hyperref}
\RequirePackage{multirow}
\RequirePackage{algorithm}
\RequirePackage{algpseudocode}
\RequirePackage{graphicx}
\RequirePackage{tcolorbox}
\RequirePackage{xcolor}
\allowdisplaybreaks

\usepackage{dpr_fec,ifthen,dsfont}

\newtheorem{assumption}{Assumption}

\newtheorem{remark}{Remark}
\newtheorem{proposition}{Proposition}

\usepackage{isetechreport}
\def\reportyear{23T}
\def\reportno{017}
\def\revisionno{0}
\def\originaldate{\today}
\coraltrue
\cvcrfalse

\begin{document}

\title {Sequential Quadratic Optimization for Stochastic Optimization with Deterministic Nonlinear Inequality and Equality Constraints}

\author{Frank E.~Curtis\thanks{E-mail: frank.e.curtis@lehigh.edu}}
\author{Daniel P.~Robinson\thanks{E-mail: daniel.p.robinson@lehigh.edu}}
\affil{Department of Industrial and Systems Engineering, Lehigh University}
\author{Baoyu Zhou\thanks{E-mail: baoyu.zhou@chicagobooth.edu}}
\affil{Booth School of Business, The University of Chicago}

\titlepage

\maketitle

%**********
% Abstract
%**********
\begin{abstract}
  A sequential quadratic optimization algorithm for minimizing an objective function defined by an expectation subject to nonlinear inequality and equality constraints is proposed, analyzed, and tested.  The context of interest is when it is tractable to evaluate constraint function and derivative values in each iteration, but it is intractable to evaluate the objective function or its derivatives in any iteration, and instead an algorithm can only make use of stochastic objective gradient estimates.  Under loose assumptions, including that the gradient estimates are unbiased, the algorithm is proved to possess convergence guarantees in expectation.  The results of numerical experiments are presented to demonstrate that the proposed algorithm can outperform an alternative approach that relies on the ability to compute more accurate gradient estimates.
\end{abstract}

%**************
% Body of paper
%**************
\newcommand{\alphatrue}{\alpha^{\rm true}}
\newcommand{\alphasuff}{\alpha^{\rm suff}}
\newcommand{\Alpha}{\Acal}
\newcommand{\Tau}{\Tcal}
\newcommand{\tautrue}{\tau^{\rm true}}
\newcommand{\tautrial}{\tau^{\rm trial}}
\newcommand{\Tautrial}{\Tcal^{\rm trial}}
\newcommand{\tautruetrial}{\tau^{\rm true,trial}}
\newcommand{\Tautruetrial}{\Tcal^{\rm true,trial}}
\newcommand{\taubartrial}{\bar\tau^{\rm trial}}
\newcommand{\tautruemin}{\tau_{\min}^{\rm true}}
\newcommand{\xitrial}{\xi^{\rm trial}}
\newcommand{\Xitrial}{\Xi^{\rm trial}}
\newcommand{\xibartrial}{\bar\xi^{\rm trial}}
\newcommand{\dtrue}{d^{\rm true}}
\newcommand{\Dtrue}{D^{\rm true}}
\newcommand{\ytrue}{y^{\rm true}}
\newcommand{\ztrue}{z^{\rm true}}
\newcommand{\lambdatrue}{\lambda^{\rm true}}
\newcommand{\utrue}{u^{\rm true}}
\renewcommand{\Null}{{\rm Null}}
\renewcommand{\argmin}{\operatornamewithlimits{argmin}}
\newcommand{\argmax}{\operatornamewithlimits{argmax}}
\newcommand{\Proj}{{\rm Proj}}

\newcommand{\fec}[1]{\textcolor{blue}{FEC: #1}}
\newcommand{\dpr}[1]{\textcolor{cyan}{#1}}
\newcommand{\bz}[1]{\textcolor{red}{#1}}

%*********
% Section
%*********
\section{Introduction}\label{sec.introduction}

We propose a sequential quadratic optimization (commonly known as SQP) algorithm for minimizing an objective function defined by an expectation subject to nonlinear inequality and equality constraints.  Such optimization problems arise in a plethora of application areas, including, but not limited to, machine learning~\cite{Lan20}, network optimization \cite{Bert98}, resource allocation \cite{Ibar88}, portfolio optimization \cite{Pero84}, risk-averse partial-differential-equation-constrained optimization~\cite{Kour16}, maximum-likelihood estimation \cite{Hath85}, and multi-stage optimization \cite{ShapDentRusz21}.

The design and analysis of deterministic algorithms for solving continuous optimization problems involving inequality and equality constraints has been a well-studied topic for decades.  Numerous types of such algorithms, such as penalty methods, interior-point methods, and SQP methods, have been designed to solve such problems.  Penalty methods are based on the idea of using unconstrained optimization algorithms to minimize a weighted sum---determined by a penalty parameter---of the objective and a measure of constraint violation; e.g., see \cite{Conn73,Flet13,Zang67} for algorithms that make use of nondifferentiable (exact) penalty functions and see \cite{DipiGrip85,DipiGrip89,GladPola79,ZavaAnit14} for algorithms that make use of differentiable (exact) penalty functions.  While they are able to offer convergence guarantees from remote starting points, the numerical performance of penalty methods often suffers from ill-conditioning of the penalty functions and/or sensitivity of the algorithm's performance on the particular scheme employed for updating the penalty parameter \cite{NoceWrig06}.  Interior-point methods~\cite{Diki67} are designed to use barrier functions to guide the algorithm along a central path through the interior of the feasible region (or, at least, the interior of a set defined by bounds on a subset of the variables) to a solution \cite{ByrdGilbNoce00,ByrdHribNoce99,LasdWareRice67,Mcgi65,Wrig97,Yama98}.  Such algorithms have been shown to be very effective in practice, which is why many state-of-the-art software packages for continuous nonlinear optimization are built on interior-point methods; see, e.g., \cite{ByrdHribNoce99,WaecBieg06}.  Overall, both penalty and interior-point methods involve the use of additional objective terms to handle the presence of inequality constraints.

Alternatively, in this paper, we present, analyze, and demonstrate the numerical performance of an SQP method for solving continuous nonlinear optimization problems.  The SQP paradigm is based on the idea of, at each iterate, solving a subproblem (or subproblems) defined based on a local linearization of the constraint function and a local quadratic approximation of the objective.  Unlike in the deterministic setting, for which numerous SQP algorithms have been proposed (see, e.g., \cite{Flet73,GillMurrSaun05,Han76,NoceWrig06}), there have been few stochastic algorithms proposed for the setting of solving optimization problems with nonlinear constraints.  That said, in the past few years, a couple of classes of stochastic SQP methods have been designed for optimization subject to nonlinear \emph{equality} constraints.  For example, the article \cite{BeraCurtRobiZhou21} proposes an SQP algorithm that uses stochastic objective gradient estimates for solving such problems that employs an adaptive step size policy based on Lipschitz constants (or estimates of them).  For an alternative setting in which one is willing to compute objective value estimates as well, and to refine objective function and gradient estimates within a given iteration until probabilistic conditions of accuracy are satisfied, the article \cite{NaAnitKola22} proposes a line-search stochastic SQP method.  There have subsequently been multiple extensions of the methods in \cite{BeraCurtRobiZhou21} and \cite{NaAnitKola22}, as well as work on different, but related algorithmic strategies---still for the setting of only nonlinear equality constraints.  There has been work on relaxing constraint qualifications~\cite{BeraCurtOneiRobi21}, allowing matrix-free and inexact solves of the arising linear systems~\cite{CurtRobiZhou21}, using a trust-region methodology \cite{FangNaMahoKola22}, incorporating noisy (potentially biased) function and gradient estimates \cite{BeraXieZhou23,OztoByrdNoce21}, employing variance-reduction strategies \cite{BeraBollZhou22,BeraShiYiZhou22}, considering sketch-and-project techniques \cite{NaMaho22}, and analyzing the worst-case complexity (see \cite{CurtOneiRobi21}) of the method proposed in \cite{BeraCurtRobiZhou21}. 

Unlike the setting of equality constraints only, to our knowledge there has been very little work on the design and analysis of stochastic algorithms for optimization subject to nonlinear (nonconvex) inequality and equality constraints.  Three exceptions are: the active-set line-search SQP algorithms proposed in~\cite{NaAnitKola21} and (very recently) in \cite{QiuKung23} and the momentum-based augmented Lagrangian method (a penalty method) proposed in \cite{ShiWangWang22}.  We expect that our proposed SQP algorithm will perform well in comparison to a stochastic-gradient-based penalty method.  We demonstrate with numerical experiments that our approach can outperform the algorithm proposed in \cite{NaAnitKola21}.  We remark in passing that interior-point methods often outperform SQP methods in the deterministic setting, but as far as we are aware there exists no interior-point method designed for the stochastic setting that we consider.

%************
% Subsection
%************
\subsection{Contributions}\label{sec.contribution}

In this paper, we build on the algorithmic strategy and analysis in \cite{BeraCurtRobiZhou21} to propose and analyze an adaptive stochastic SQP algorithm for solving nonlinear optimization problems subject to (deterministic) inequality and equality constraints.  This work involves significant advancements beyond that in \cite{BeraCurtRobiZhou21} that are necessary since, unlike in the setting of only having equality constraints, the presence of inequality constraints automatically guarantees that, at a given iterate, the search direction computed in a stochastic SQP method will be a biased estimate of the ``true'' search direction, i.e., the one that would be computed if the actual gradient of the objective function were available.  This necessitates a distinct change in the design of the algorithm as well as distinct alterations to the convergence analysis, since the analysis in \cite{BeraCurtRobiZhou21} relies heavily on the search directions being (conditionally) unbiased estimators of their ``true'' counterparts.  The algorithm from the literature that can be seen as the nearest alternative approach is the algorithm in~\cite{NaAnitKola21}.  However, there are substantial differences between the algorithm and analysis in \cite{NaAnitKola21} and those presented in this paper.  Like in \cite{NaAnitKola22} for the equality-only case, the algorithm in~\cite{NaAnitKola21} is designed for the setting in which one is willing to refine function and gradient estimates within an iteration until probabilistic conditions of accuracy are satisfied, and in this manner the analysis of that algorithm offers guarantees that are relatively closer to those offered for a deterministic algorithm.  By contrast, the algorithm in this paper, like the algorithm in \cite{BeraCurtRobiZhou21}, is designed to allow the stochastic gradient estimates to be potentially much less accurate, and in such a context we are satisfied with offering convergence guarantees in expectation.  We compare the numerical performance of our proposed algorithm with that in \cite{NaAnitKola21} to demonstrate that there are settings in which our proposed approach has advantages in practice.

%************
% Subsection
%************
\subsection{Notation}\label{sec.notation}

We use $\R{}$ to denote the set of real numbers, $\Rext{}$ to denote the set of extended-real numbers (i.e., $\Rext{} := \R{} \cup \{-\infty,\infty\}$), and $\R{}_{\geq a}$ (resp.,~$\R{}_{>a}$) to denote the set of real numbers greater than or equal to (resp.,~greater than) $a \in \R{}$.  We append a superscript to such a set to denote the space of vectors or matrices whose elements are restricted to the indicated set; e.g.,~we use $\R{n}$ to denote the set of $n$-dimensional real vectors and~$\R{m \times n}$ to denote the set of $m$-by-$n$-dimensional real matrices.  We use $\N{} := \{1,2,\dots\}$ to denote the set of positive integers and, given $n \in \N{}$, we use~$[n] := \{1,\dots,n\}$ to denote the set of positive integers less than or equal to $n$.  Given $(a,b) \in \R{n} \times \R{n}$, we write $a \perp b$ to mean, with $a_i$ and $b_i$ denoting the $i$th elements of $a$ and $b$, respectively, that $a_i = 0$ and/or $b_i = 0$ for all $i \in [n]$.  Given real symmetric matrices $A \in \R{n \times n}$ and $B \in \R{n \times n}$, we write $A \succeq B$ (resp.,~$A \succ B$) to indicate that $A - B$ is positive semidefinite (resp., positive definite).  Given $H \in \R{n \times n}$ with $H \succ 0$ and $a \in \R{n}$, we denote the norm $\|a\|_H := \sqrt{a^THa}$.

Our problem of interest is defined with respect to a variable $x \in \R{n}$ and the algorithm that we propose and analyze is iterative, meaning that, in any run, it generates an iterate sequence that we denote as $\{x_k\}$ with $x_k \in \R{n}$ for all generated $k \in \N{}$, i.e., $\{x_k\} \subset \R{n}$.  (We use such notation throughout the paper when the elements of sequence are contained within a given set.  We say ``for all \emph{generated} $k \in \N{}$'' since our proposed algorithm might terminate finitely.  Whether a subscript is being used to indicate the element of a vector or the index number of a sequence is always made clear by the context.  The $i$th element of an iterate $x_k$ is denoted $[x_k]_i$.)  We use subscripts similarly to denote other quantities corresponding to each iteration of the algorithm; e.g., we introduce a merit parameter denoted as $\tau \in \R{}_{>0}$ whose value in iteration $k \in \N{}$ is denoted as $\tau_k \in \R{}_{>0}$, and corresponding to a constraint function $c$ (see problem~\eqref{prob.opt} below) we denote its value at $x_k$ as $c_k := c(x_k)$.

The iteration-dependent quantities mentioned in the previous paragraph---and additional ones introduced in the description of our algorithm---represent realizations of the random variables in a stochastic process generated by the algorithm.  Specifically, the behavior of our algorithm is dictated by prescribed initial conditions and a sequence of stochastic objective gradient estimators that we denote by $\{G_k\}$.  After proving preliminary results that hold for every run of the algorithm, we present our ultimate convergence theory for our algorithm in terms of a filtration defined in terms of $\sigma$-algebras dependent on the initial conditions of the algorithm and $\{G_k\}$.

%************
% Subsection
%************
\subsection{Organization}\label{sec.organization}

A statement of our problem of interest and preliminary assumptions about its objective and constraint functions, as well as about user-defined quantities in our proposed algorithm, are stated in Section~\ref{sec.setting}.  A description of our proposed algorithm is provided in Section~\ref{sec.algorithm}.  Convergence-in-expectation of the algorithm is proved under reasonable assumptions in Section~\ref{sec.analysis}.  The results of numerical experiments are presented in Section~\ref{sec.numerical} and concluding remarks are given in Section~\ref{sec.conclusion}.

%*********
% Section
%*********
\section{Setting}\label{sec.setting}

We formulate our problem of interest as
\bequation\label{prob.opt}
  \min_{x\in\R{n}}\ f(x)\ \text{subject to (s.t.)}\ c(x) = 0\ \text{and}\ x \geq 0\ \text{with}\ f(x) = \E_\omega[F(x,\omega)],
\eequation
where $f : \R{n} \to \R{}$ and $c : \R{n} \to \R{m}$ are continuously differentiable, $\omega$ is a random variable with associated probability space $(\Omega,\Fcal,\P)$, $F : \R{n} \times \Omega \to \R{}$, and $\E_\omega[\cdot]$ denotes expectation taken with respect to the distribution of $\omega$.  Our algorithm and analysis extend easily to the setting in which the nonnegativity constraint in \eqref{prob.opt} is generalized to $l \leq x \leq u$ for some $(l,u) \in \Rext{n} \times \Rext{n}$ with $l_i \leq u_i$ for all $i \in [n]$; we merely consider nonnegativity in \eqref{prob.opt} for the sake of notational simplicity.  It is also worth mentioning that any smooth constrained optimization problem can be reformulated as \eqref{prob.opt} (or at least as such a problem with generalized bound constraints); e.g., inequality constraints $c_\Ical(x) \leq 0$, where $c_{\Ical} : \R{n} \to \R{m_{\Ical}}$ is continuously differentiable, can be reformulated to fit into the form of \eqref{prob.opt} through the incorporation of slack variables, say $s \in \R{m_{\Ical}}$, to have the constraints $c_{\Ical}(x) + s_{\Ical} = 0$ and $s_{\Ical} \geq 0$.

We make the following assumption throughout the remainder of the paper pertaining to the functions in problem~\eqref{prob.opt} and our proposed algorithm.  As seen in the following section, our algorithm seeks feasibility and stationarity with respect to~\eqref{prob.opt} by generating an iterate sequence that stays feasible with respect to the bound constraints, meaning that, in any run of the algorithm, $x_k \in \R{n}_{\geq0}$ for all generated $k \in \N{}$.

\begin{assumption}\label{ass.prob}
  Let $\Xcal \subset \R{n}$ be an open convex set that almost-surely contains the iterate sequence $\{x_k\} \subset \R{n}_{\geq 0}$ generated in any realization of a run of the algorithm.  The objective function $f : \R{n} \to \R{}$ is continuously differentiable and bounded below over $\Xcal$ and the objective gradient function $\nabla f : \R{n} \to \R{n}$ is Lipschitz continuous and bounded in norm over~$\Xcal$.  Similarly, for all $i \in [m]$, the constraint function $c_i : \R{n} \to \R{}$ is continuously differentiable and bounded over $\Xcal$ and the constraint gradient function $\nabla c_i : \R{n} \to \R{n}$ is Lipschitz continuous and bounded in norm over $\Xcal$.  Finally, the constraint Jacobian $\nabla c^T : \R{n} \to \R{m \times n}$ has full row rank over $\Xcal$.
\end{assumption}

Under Assumption~\ref{ass.prob}, there exists $f_{\inf} \in \R{}$ and a tuple of positive constants $(\kappa_{\nabla f}, \kappa_c, \kappa_{\nabla c}, L, \Gamma) \in \R{}_{>0} \times \R{}_{>0} \times \R{}_{>0} \times \R{}_{>0} \times \R{}_{>0}$ such that for all $x \in \Xcal$ one has
\bequation\label{eq.bounds}
  f(x) \geq f_{\inf},\ \ \|\nabla f(x)\|_2 \leq \kappa_{\nabla f},\ \ \|c(x)\|_2 \leq \kappa_c,\ \ \text{and}\ \ \|\nabla c(x)\|_2 \leq \kappa_{\nabla c},
\eequation
and for all $(x,\xbar) \in \Xcal \times \Xcal$ one has
\bequation\label{eq.Lipschitz}
  \|\nabla f(x) - \nabla f(\xbar)\|_2 \leq L \|x - \xbar\|_2\ \ \text{and}\ \ \|\nabla c(x)^T - \nabla c(\xbar)^T\|_2 \leq \Gamma \|x - \xbar\|_2.
\eequation
In addition, due to the continuous differentiability of the objective and constraint functions and the full row rank of the constraint Jacobians, it follows that at any (local) minimizer of \eqref{prob.opt}, call it $x \in \R{n}$, there exists $y \in \R{m}$ and $z \in \R{n}$ such that the following Karush-Kuhn-Tucker (KKT) conditions are satisfied:
\bequation\label{eq.KKT}
  \nabla f(x) + \nabla c(x)y - z = 0,\ \ c(x) = 0,\ \ 0 \leq x \perp z \geq 0.
\eequation
We refer to any $x \in \R{n}$ such that there exists $(y,z) \in \R{m} \times \R{n}$ satisfying \eqref{eq.KKT} as a first-order stationary point (or KKT point) with respect to \eqref{prob.opt}.

Since our algorithm generates iterates that are feasible with respect to the bound constraints, but not necessarily with respect to the equality constraints, we need to account for the possible existence of points that are infeasible for \eqref{prob.opt}, but are stationary with respect to the minimization of a constraint violation measure over~$\R{n}_{\geq0}$.  We refer to a point that is infeasible for \eqref{prob.opt} as an infeasible stationary point if it is stationary with respect to the minimization of $\thalf \|c(x)\|_2^2$ subject to $x \in \R{n}_{\geq0}$, meaning
\bequation\label{eq.KKT_con}
  0 \leq x \perp \nabla c(x)c(x) \geq 0.
\eequation

Each iteration of our algorithm requires a stochastic estimate of the gradient of the objective at the current iterate.  In a given run at iteration $k \in \N{}$, the realization of the iterate and gradient estimate is $(x_k,g_k)$, which later in our analysis we denote as a realization of the pair of random variables $(X_k,G_k)$.  (See Section~\ref{sec.guarantees} for a complete description of a stochastic process that we analyze.)  With respect to the gradient estimators, we make Assumption~\ref{ass.g} below.  For the prescribed (i.e., not random) sequence $\{\rho_k\} \subset \R{}_{>0}$ referenced in the assumption, we state precise conditions that it must satisfy in Section~\ref{sec.guarantees}.  In the assumption and throughout the remainder of the paper, we use $\E_k[\cdot]$ to denote expectation taken with respect to the distribution of $\omega$ conditioned on a trace $\sigma$-algebra of an event $\Ecal$, denoted by $\Fcal_k$; see Section~\ref{sec.guarantees}.

\bassumption\label{ass.g}
  For a prescribed $\{\rho_k\} \subset \R{}_{>0}$, one finds for all $k \in \N{}$ that
  \bequation\label{eq.g}
    \E_k[G_k] = \nabla f(X_k)\ \ \text{and}\ \ \E_k[\|G_k - \nabla f(X_k)\|_2^2] \leq \rho_k.
  \eequation
\eassumption
One might relax the latter condition in \eqref{eq.g} and obtain guarantees that are similar to those that we prove; see, e.g., \cite{PateZhan21}. We employ \eqref{eq.g} for simplicity, since it is sufficient for demonstrating the guarantees that our algorithmic approach can offer.

Each iteration of our algorithm also makes use of a symmetric and positive-definite (SPD) matrix, denoted as $H_k \in \R{n \times n}$ for iteration $k \in \N{}$, to define a quadratic term in the subproblem that is solved for computing the search direction.  For simplicity, we assume that the sequence $\{H_k\}$ is prescribed, e.g., one may consider $H_k = I$ for all $k \in \N{}$.  More generally, one could consider a more sophisticated scheme such as setting, for all $k \in \N{}$, the matrix $H_k$ as a stochastic estimate of the Hessian of the objective function and/or a Lagrangian function as long as it is sufficiently positive definite and bounded and the choice is made to be conditionally uncorrelated with the stochastic gradient estimate.  However, since considering such a loose requirement would only obfuscate our analysis without adding significant value, we assume for simplicity that $\{H_k\}$ is prescribed and merely satisfies the following.

\bassumption\label{ass.H}
  There exists $(\kappa_H,\zeta) \in \R{}_{>0} \times \R{}_{>0}$ with $\kappa_H \geq \zeta$ such that, for all $k \in \N{}$, the SPD matrix $H_k \in \R{n \times n}$ has $\kappa_H I \succeq H_k \succeq \zeta I$.
\eassumption

Observe from Assumption~\ref{ass.H} that we are not assuming that accurate second-order information is being used by the algorithm.  Hence, our convergence guarantees are of the type that may be expected for a first-order-type algorithm, although in situations when it is computationally tractable, one might find better performance if~$H_k$ incorporates some (approximate) second-order derivative information.

%*********
% Section
%*********
\section{Algorithm}\label{sec.algorithm}

In this section, we present our proposed algorithm.  We state the algorithm in terms of a particular realization of it (e.g., denoting the iterate for each $k \in \N{}$ as $x_k$), although our subsequent analysis of it (starting in Section~\ref{sec.guarantees}) will be written in terms of the stochastic process that the algorithm defines.

Each iteration of our algorithm proceeds as follows.  First, given the current iterate $x_k \in \R{n}_{\geq0}$, the algorithm computes a direction whose purpose is to determine the progress that can be made in terms of reducing a measure of violation of a linearization of the equality constraints subject to the bound constraints.  This is done in a manner that regularizes the component of the direction that lies in the null space of the constraint Jacobian.  Specifically, the iteration commences by computing a direction $v_k := u_k + \nabla c(x_k) w_k \in \R{n}$, where $u_k \in \Null(\nabla c(x_k)^T)$ and $\nabla c(x_k) w_k \in \Range(\nabla c(x_k))$, by solving the quadratic optimization subproblem
\bequation\label{prob.v}
  \baligned
    \min_{u \in \R{n}, w \in \R{m}} &\ \thalf \|c_k + \nabla c(x_k)^T\nabla c(x_k)w \|_2^2 + \thalf \mu_k \|u\|_2^2 \\
    \st &\ \nabla c(x_k)^Tu = 0\ \ \text{and}\ \ x_k + u + \nabla c(x_k)w \geq 0,
  \ealigned
\eequation
where $\mu_k \in \R{}_{>0}$ is a user-prescribed parameter.  Observe that since $x_k \in \R{n}_{\geq0}$, this subproblem is always feasible, and by construction it is convex.  Generally, the solution of \eqref{prob.v} might not be unique, but in our setting it is unique since $\nabla c(x_k)^T$ has full row rank.  In our analysis, we show that the solution of subproblem~\eqref{prob.v} is given by $(u_k,w_k) = (0,0)$ if and only if the current iterate $x_k$ is stationary with respect to the minimization of $\thalf \|c(x)\|_2^2$ over $x \in \R{n}_{\geq0}$.  This means, e.g., that if $c_k \neq 0$, but the solution of~\eqref{prob.v} is $(u_k,w_k) = (0,0)$---which, by the Fundamental Theorem of Linear Algebra, occurs if and only if $v_k = u_k + \nabla c(x_k) w_k = 0$---then it is reasonable to terminate since $x_k$ is an infeasible stationary point (see \eqref{eq.KKT_con}), as in our algorithm.

After computing $v_k \in \R{n}$ by solving \eqref{prob.v} and generating a stochastic objective gradient estimate~$g_k \in \R{n}$ (see Assumption~\ref{ass.g}), the algorithm next computes a search direction $d_k \in \R{n}$ by solving the quadratic optimization subproblem
\bequation\label{prob.d}
  \min_{d \in \R{n}}\ g_k^Td + \thalf d^TH_kd\ \st\ \nabla c(x_k)^Td = \nabla c(x_k)^Tv_k\ \text{and}\ x_k + d \geq 0.
\eequation
By construction, this subproblem is feasible, and under Assumption~\ref{ass.H} it is convex.  The search direction $d_k$ is designed to achieve the same progress toward linearized feasibility within the nonnegative orthant that is achieved by $v_k$, then within the null space of $\nabla c(x_k)^T$ and the nonnegative orthant aims to minimize a (stochastically estimated) local quadratic approximation of the objective function at $x_k$.

The remainder of the $k$th iteration proceeds in a similar manner as in \cite{BeraCurtRobiZhou21,CurtRobiZhou21}.  In particular, with the $\ell_2$-norm merit function in mind, namely, $\phi : \R{n} \times \R{}_{>0} \to\R{}$ defined by $\phi(x,\tau) = \tau f(x) + \|c(x)\|_2$, the algorithm next sets a value for the merit parameter $\tau_k \in \R{}_{>0}$.  This is done by considering a local model of this merit function, namely, $l : \R{n} \times \R{}_{>0} \times \R{n} \times \R{n} \to \R{}$ defined by $l(x,\tau,g,d) = \tau(f(x) + g^Td) + \|c(x) + \nabla c(x)^Td\|_2$, and in particular the reduction in this model defined for all $k \in \N{}$ by
\bequation\label{eq.model_reduction}
  \baligned
    \Delta l(x_k,\tau_k,g_k,d_k)
    :=&\ l(x_k,\tau_k,g_k,0) - l(x_k,\tau_k,g_k,d_k) \\
    =&\ -\tau_k g_k^Td_k + \|c_k\|_2 - \|c_k + \nabla c(x_k)^Td_k\|_2,
  \ealigned
\eequation
and setting $\tau_k$ such that this reduction is sufficiently large.  Specifically, if $d_k \neq 0$, then with user-prescribed $(\epsilon_\tau,\sigma) \in (0,1) \times (0,1)$, the algorithm first sets
\bequation\label{eq.tautrial}
  \tautrial_k \gets \bcases \infty & \text{if $g_k^Td_k + \thalf d_k^TH_kd_k \leq 0$} \\ \frac{(1-\sigma)(\|c_k\|_2 - \|c_k + \nabla c(x_k)^Td_k\|_2)}{g_k^Td_k + \thalf d_k^TH_kd_k} & \text{otherwise,} \ecases 
\eequation
then sets the merit parameter value as
\bequation\label{eq.merit_parameter}
  \tau_k \gets \bcases \tau_{k-1} & \text{if $\tau_{k-1} \leq \tautrial_k$} \\ \min\{(1-\epsilon_{\tau}) \tau_{k-1}, \tautrial_k\} & \text{otherwise.} \ecases   
\eequation
(The value $\tau_0 \in \R{}_{>0}$ is also prescribed by the user.)  On the other hand, if $d_k = 0$, then the algorithm simply sets $\tautrial_k \gets \infty$ and $\tau_k \gets \tau_{k-1}$.  We show in our analysis (see Lemma~\ref{lem.model_reduction}) that this procedure for setting $\tau_k$ ensures that $\Delta l(x_k, \tau_k, g_k, d_k)$ is sufficiently large relative to the squared norm of the search direction and the improvement offered toward linearized feasibility.  For use in the step size procedure, the algorithm next sets a value $\xi_k \in \R{}_{>0}$ (referred to as the ratio parameter) that acts as an estimate for a lower bound of the ratio between the model reduction and a multiple of the squared norm of the search direction.  Specifically, if $d_k \neq 0$, it sets
\bequation\label{eq.xi_trial}
  \xitrial_k \gets \tfrac{\Delta l(x_k,\tau_k,g_k,d_k)}{\tau_k\|d_k\|_2^2},\ \ \text{then}\ \ \xi_k \gets \bcases \xi_{k-1} & \text{if $\xi_{k-1} \leq \xitrial_k$} \\ \min\{(1-\epsilon_{\xi})\xi_{k-1},\xitrial_k\} & \text{otherwise,} \ecases
\eequation
where $(\xi_0,\epsilon_{\xi}) \in \R{}_{>0} \times (0,1)$ are user-prescribed parameters; see \cite{BeraCurtRobiZhou21,CurtRobiZhou21} for further motivation.  On the other hand, if $d_k = 0$, then it sets $\xitrial_k \gets \infty$ and $\xi_k \gets \xi_{k-1}$.

The step size selection procedure, which for all $k \in \N{}$ chooses the step size $\alpha_k \in \R{}_{>0}$, can now be summarized as follows.  First, suppose that $d_k \neq 0$. With user-prescribed $\eta\in(0,1)$, $\theta \in \R{}_{>0}$, and $\{\beta_k\}$ with $\beta_k \in (0,1]$ for all $k \in \N{}$ such that
\bequation\label{eq.beta_constraint}
  \alpha_k^{\min} \gets \tfrac{2 (1-\eta) \beta_k \xi_k \tau_k}{\tau_k L + \Gamma} \in (0,1]\ \ \text{for all}\ \ k\in\N{},
\eequation
and with the strongly convex function $\varphi_k : \R{}_{\geq0} \to \R{}$ defined by
\bequation\label{eq.varphi}
  \baligned
    \varphi_k(\alpha)
    =&\ (\eta - 1) \alpha \beta_k \Delta l(x_k,\tau_k,g_k,d_k) + \|c_k + \alpha \nabla c(x_k)^T d_k\|_2 - \|c_k\|_2 \\
    &\ + \alpha(\|c_k\|_2 - \|c_k + \nabla c(x_k)^T d_k\|_2) + \thalf(\tau_kL + \Gamma)\alpha^2\|d_k\|_2^2,
  \ealigned
\eequation
the algorithm sets the values
\bequation\label{eq.alpha_max}
  \baligned
    \alpha_k^{\varphi} &\gets \max\{\alpha\in\R{}_{\geq 0}:\varphi_k(\alpha) \leq 0\} \\
    \text{and}\ \ \alpha_k^{\max} &\gets \min \{1,\alpha_k^{\varphi},\alpha_k^{\min} + \theta\beta_k \}.
  \ealigned
\eequation
The algorithm then chooses the step size $\alpha_k$ as any value in $[\alpha_k^{\min},\alpha_k^{\max}]$.  Second, if $d_k = 0$, then the algorithm simply sets all step size values to 1.

A complete statement of our algorithm is given as Algorithm~\ref{alg.sqp}. 

\begin{algorithm}[ht]
  \caption{Stochastic SQP}
  \label{alg.sqp}
  \begin{algorithmic}[1]
    \Require $x_1\in\R{n}_{\geq0}$; $\{\mu_k\} \subset \R{}_{>0}$; $\{H_k\} \subset \R{n \times n}$ satisfying Assumption~\ref{ass.H}; $\tau_0 \in \R{}_{>0}$; $\xi_0 \in \R{}_{>0}$; $\{\sigma,\eta,\epsilon_\tau,\epsilon_{\xi}\} \subset (0,1)$; $\{\beta_k\} \subset (0,1]$ satisfying \eqref{eq.beta_constraint}; $\theta \in \R{}_{>0}$; $\{\rho_k\} \subset \R{}_{>0}$; Lipschitz constants $L \in \R{}_{>0}$ and $\Gamma \in \R{}_{>0}$ (see \eqref{eq.Lipschitz}) 
    \For{$k\in\N{}$}
      \State compute $v_k \in \R{n}$ by solving \eqref{prob.v}
      \State \textbf{if} $c_k \neq 0$ and $v_k = 0$ \textbf{then terminate} and \textbf{return} $x_k$ (infeasible stationary)
      \State compute $g_k \in \R{n}$ (recall Assumption~\ref{ass.g}) \label{step.g}
      \State compute $d_k \in \R{n}$ by solving \eqref{prob.d}
      \If{$d_k = 0$} 
        \State set $\tautrial_k \gets \infty$, $\tau_k \gets \tau_{k-1}$, $\xitrial_k \gets \infty$, and $\xi_k \gets \xi_{k-1}$
        \State set $\alpha_k^{\min} \gets 1$, $\alpha_k^\varphi \gets 1$, $\alpha_k^{\max} \gets 1$, and $\alpha_k \gets 1$
      \Else
        \State set $\tautrial_k$ by \eqref{eq.tautrial} and $\tau_k$ by \eqref{eq.merit_parameter}
        \State set $\xitrial_k$ and $\xi_k$ by \eqref{eq.xi_trial}
        \State set $\alpha_k^{\min}$ by \eqref{eq.beta_constraint} and both $\alpha_k^\varphi$ and $\alpha_k^{\max}$ by \eqref{eq.alpha_max}
        \State choose $\alpha_k\in[\alpha_k^{\min},\alpha_k^{\max}]$
      \EndIf
    \State Set $x_{k+1}\gets x_k + \alpha_k d_k$ 
    \EndFor
  \end{algorithmic}
\end{algorithm}

%*********
% Section
%*********
\section{Analysis}\label{sec.analysis}

In this section, we provide theoretical results for Algorithm~\ref{alg.sqp}.  We begin by introducing common assumptions under which one can establish stationarity measures for problem~\eqref{prob.opt} that are defined by solutions of~\eqref{prob.v} and/or \eqref{prob.d}.  These stationarity measures allow us to connect our convergence guarantees for Algorithm~\ref{alg.sqp} with stationarity conditions for \eqref{prob.opt}.  Then, under Assumptions~\ref{ass.prob} and~\ref{ass.H}, we prove generally applicable results pertaining to the behavior of algorithmic quantities in any run of the algorithm.  These results reveal that the algorithm is well defined in the sense that any run will either terminate and return an infeasible stationary point or generate an infinite sequence of iterates.  We then consider convergence properties of the algorithm in the event that the (monotonically nonincreasing) merit parameter sequence eventually produces values that are sufficiently small, yet bounded away from zero, which, as shown in our analysis, means that the sequence ultimately becomes constant at a sufficiently small value.  This analysis, which includes our main convergence results for the algorithm, is provided under Assumption~\ref{ass.good_merit_parameter} stated on page~\pageref{ass.good_merit_parameter}.  We follow this analysis with a section on theoretical results related to the occurrence of the event in Assumption~\ref{ass.good_merit_parameter}.  As in~\cite{BeraCurtRobiZhou21} for the equality-constraints-only setting, this discussion illuminates the fact that while the event in Assumption~\ref{ass.good_merit_parameter} is not always guaranteed to occur due to the looseness of our assumptions about properties of the stochastic gradient estimates, the event represents likely behavior in practice, which shows that our convergence results about the algorithm are meaningful for real-world situations.  We conclude this section with a discussion of the behavior of the algorithm in the deterministic setting, i.e., when the true gradient of the objective is employed in all iterations.  This discussion is meant to provide confidence to a user that our algorithm is based on one that has state-of-the-art convergence properties under common assumptions in the deterministic setting.

%************
% Subsection
%************
\subsection{Subproblems and Stationarity Measures}\label{sec.subproblems}

We begin by showing that subproblem~\eqref{prob.v} yields a zero solution if and only if the point defining the subproblem is feasible for problem~\eqref{prob.opt} or an infeasible stationary point.

\blemma\label{lem.normal_generic}
  Suppose that Assumption~\ref{ass.prob} holds, $x \in \Xcal \cap \R{n}_{\geq0}$, and, given $\mu \in \R{}_{>0}$, consider the quadratic optimization problem $($recall \eqref{prob.v}$)$
  \bequation\label{prob.v_generic}
    \baligned
      \min_{u\in\R{n},w\in\R{m}} &\ \thalf \|c(x) + \nabla c(x)^T\nabla c(x)w\|_2^2 + \thalf \mu \|u\|_2^2 \\
      \st &\ \nabla c(x)^T u = 0\ \ \text{and}\ \ x + u + \nabla c(x)w \geq 0.
    \ealigned
  \eequation
  Then, the unique optimal solution of problem~\eqref{prob.v_generic} is $(u,w) = (0,0)$ if and only if $x$ is feasible for problem~\eqref{prob.opt} or an infeasible stationary point $($i.e., it satisfies \eqref{eq.KKT_con}$)$, whereas $(u,w) \neq (0,0)$ if and only if $\|c(x)\|_2 > \|c(x) + \nabla c(x)^T \nabla c(x) w\|_2$.
\elemma
\bproof
  Suppose the conditions of the lemma hold and let $(u,w)$ be the unique optimal solution of \eqref{prob.v_generic}. Since $x \in \R{n}_{\geq0}$, it follows that $(0,0)$ is feasible for~\eqref{prob.v_generic}.  In addition, necessary and sufficient optimality conditions for~\eqref{prob.v_generic} are that, corresponding to $(u,w) \in \R{n} \times \R{m}$, there exists $(\gamma,\delta) \in \R{m} \times \R{n}$ with
  \bequation\label{eq.KKT_con_sub_generic}
    \baligned
      \nabla c(x)^T \nabla c(x) c(x) + \nabla c(x)^T \nabla c(x) \nabla c(x)^T \nabla c(x) w - \nabla c(x)^T\delta &= 0, \\
      \mu u + \nabla c(x) \gamma - \delta = 0,\ \ \nabla c(x)^T u = 0,\ \ \text{and}\ \ 
      0 \leq \delta \perp x + u + \nabla c(x) w &\geq 0.
    \ealigned
  \eequation
  If $(u,w) = (0,0)$, then it follows from \eqref{eq.KKT_con_sub_generic} that
  \bequation\label{eq.KKT_con_sub_zero_generic}
    \nabla c(x)^T \nabla c(x) c(x) - \nabla c(x)^T \delta = 0,\ \ \nabla c(x) \gamma - \delta = 0,\ \ \text{and}\ \ 0 \leq \delta \perp x \geq 0.
  \eequation
  Since $\nabla c(x)^T$ has full row rank, \eqref{eq.KKT_con_sub_zero_generic} implies $\gamma = (\nabla c(x)^T \nabla c(x))^{-1} \nabla c(x)^T \delta = c(x)$, $\delta = \nabla c(x) c(x)$, and $0 \leq \nabla c(x) c(x) \perp x \geq 0$, which from \eqref{eq.KKT_con} means that $x$ is either feasible or an infeasible stationary point, as desired.  On the other hand, if $x$ is either feasible or an infeasible stationary point, meaning $0 \leq \nabla c(x) c(x) \perp x \geq 0$, then $u = 0$, $w = 0$, $\gamma = c(x)$, and $\delta = \nabla c(x) c(x)$ satisfy \eqref{eq.KKT_con_sub_generic}, and this solution (i.e., $(u,w) = (0,0)$) is unique since the objective of~\eqref{prob.v_generic} is strongly convex.

  Now let us show that the unique optimal solution of~\eqref{prob.v_generic} is $(u,w) \neq (0,0)$ if and only if $\|c(x)\|_2 > \|c(x) + \nabla c(x)^T \nabla c(x) w\|_2$.  If $\|c(x)\|_2 > \|c(x) + \nabla c(x)^T \nabla c(x) w\|_2$, then $w \neq 0$ follows trivially, giving the desired conclusion.  To prove the reverse implication, let us consider two cases.  If $u \neq 0$, then, since~$(0,0)$ is feasible for~\eqref{prob.v_generic},
  \bequationNN
     \thalf \|c(x)\|_2^2 \geq \thalf \|c(x) + \nabla c(x)^T \nabla c(x) w\|_2^2 + \thalf \mu \|u\|_2^2 > \thalf \|c(x) + \nabla c(x)^T \nabla c(x) w\|_2^2,
  \eequationNN
  as desired.  Second, if $u = 0$ and $w \neq 0$, then $w$ is the minimizer of the strongly convex objective $\thalf \|c(x) + \nabla c(x)^T \nabla c(x) w\|_2^2$ subject to $x + \nabla c(x) w \geq 0$.  Since $0$ is feasible for this problem, $w \neq 0$ means that $\thalf \|c(x)\|_2^2 > \thalf \|c(x) + \nabla c(x)^T \nabla c(x) w\|_2^2$, as desired.
\eproof

We now show that, under common assumptions and given $x_k \in \R{n}_{\geq0}$, the quantity $\|v_k\|_2^2$, where $v_k \in \R{n}$ solves subproblem \eqref{prob.v}, represents a stationarity measure with respect to the problem to minimize $\thalf \|c(x)\|_2^2$ subject to $x \in \R{n}_{\geq0}$.  (The assumption in the lemma that $\mu_k = \mu \in \R{}_{>0}$ for all $k \in \N{}$ could be relaxed; see Remark~\ref{rem.mu_H} at the end of this subsection.  We consider this case for the sake of brevity.)

\blemma\label{lem.normal_stationarity}
  Suppose that Assumption~\ref{ass.prob} holds and there exists infinite $\Scal \subseteq \N{}$ such that for some sequence $\{x_k\} \subset \Xcal \cap \R{n}_{\geq0}$ one finds $\{x_k\}_{k \in \Scal} \to x_*$ for some $x_* \in \Xcal \cap \R{n}_{\geq0}$ where, with $\Acal(x) := \{i \in [n] : x_i = 0\}$, $I_{\Acal(x)}$ denoting the matrix composed of rows of $I \in \R{n \times n}$ corresponding to indices in $\Acal(x)$, and $\nabla c(x)_{\Acal(x)}$ denoting the matrix composed of rows of $\nabla c(x)$ corresponding to indices in $\Acal(x)$, one finds that
  \benumerate
    \item[(i)] $[\nabla c(x_*) c(x_*)]_i > 0$ for all $i \in \Acal(x_*)$ and
    \item[(ii)] the following matrix has full row rank:
    %\bequationNN
      $\bbmatrix 0 & \nabla c(x_*)^T \\ \nabla c(x_*)_{\Acal(x_*)} & I_{\Acal(x_*)} \ebmatrix$.
    %\eequationNN
  \eenumerate
  Then, with $\mu_k = \mu \in \R{}_{>0}$ for all $k \in \N{}$, and with $(u_k,w_k)$ solving subproblem~\eqref{prob.v} and $v_k := u_k + \nabla c(x) w_k$ for all $k \in \N{}$, it follows that $x_*$ satisfies the stationarity conditions \eqref{eq.KKT_con} if and only if $\{v_k\}_{k \in \Scal} \to 0$.
\elemma
\bproof
  Let $\Acal_* := \Acal(x_*)$ and $j(x) := \nabla c(x)^T$ and consider the linear system
  \bequationNN
    \bbmatrix j(x) j(x)^T j(x) j(x)^T & 0 & 0 & -j(x)_{\Acal_*} \\ 0 & \mu I & j(x)^T & -I_{\Acal_*}^T \\ 0 & j(x) & 0 & 0 \\ j(x)_{\Acal_*}^T & I_{\Acal_*} & 0 & 0 \ebmatrix \bbmatrix w \\ u \\ \gamma \\ \delta_{\Acal_*} \ebmatrix = \bbmatrix -j(x) j(x)^T c(x) \\ 0 \\ 0 \\ -x_{\Acal_*} \ebmatrix.
  \eequationNN
  Since, under the conditions of the lemma, the matrix in this linear system is nonsingular when $x = x_*$ (e.g., this follows from~\cite[Theorem 1.5.1]{Robi07} and (ii)), it follows that there exists an open ball $\Bcal_*$ centered at $x_*$ such that, for each $x \in \Bcal_* \cap \Xcal \cap \R{n}_{\geq0}$, this linear system has a unique solution, call it $(w(x),u(x),\gamma(x),\delta_{\Acal_*}(x))$, and---due to continuity of the left-hand-side matrix and right-hand-side vector with respect to~$x$---this solution varies continuously over $\Bcal_* \cap \Xcal \cap \R{n}_{\geq0}$.  If $x_*$ satisfies~\eqref{eq.KKT_con}, then it follows that $(0,0,c(x_*),[j(x_*)^T c(x_*)]_{\Acal_*})$ (with $[j(x_*)^T c(x_*)]_{\Acal_*} > 0$) is the unique solution of the system at $x = x_*$, and for all $x \in \Bcal_* \cap \Xcal \cap \R{n}_{\geq0}$ the solution of the system in conjunction with $\delta_i = 0$ for all $i \notin \Acal_*$ satisfies \eqref{eq.KKT_con_sub_generic}, meaning that the components $(u(x),w(x))$ represent the unique optimal solution of problem~\eqref{prob.v_generic}.  Hence, with respect to the quantities in the lemma and using Assumption~\ref{ass.prob}, one finds that $\{v_k\}_{k \in \Scal} \to 0$, as desired. To prove the reverse inclusion, suppose that $\{v_k\}_{k \in \Scal} \to 0$, from which it follows by the Fundamental Theorem of Linear Algebra and (ii) that $\{(u_k,w_k)\}_{k \in \Scal} \to 0$.  For all $k \in \Scal$, let $(u_k,w_k,\gamma_k,\delta_k)$ be a primal-dual optimal solution of \eqref{prob.v} (satisfying optimality conditions of the form in \eqref{eq.KKT_con_sub_generic}).  One finds under the conditions of the lemma that, for all sufficiently large $k \in \Scal$, this solution has $[\delta_k]_i = 0$ for all $i \notin \Acal(x_*)$ whereas $(u_k,w_k,\gamma_k,[\delta_k]_{\Acal_*})$ solves the linear system above at $x = x_k$.  Since, by the arguments above, this solution varies continuously within $\Bcal_* \cap \Xcal \cap \R{n}_{\geq0}$, the fact that $\{x_k\}_{k \in \Scal} \to x_*$ implies that $x_*$ satisfies~\eqref{eq.KKT_con}, as desired.
\eproof

In fact, under the conditions of the prior lemma, the quantity $\|c_k\|_2 - \|c_k + \nabla c(x_k)^T v_k\|_2$ also represents a stationarity measure for the problem to minimize $\thalf \|c(x)\|_2^2$ subject to $x \in \R{n}_{\geq0}$.  This is shown in the following lemma.

\blemma\label{lem.normal_stationarity_2}
  Suppose that Assumption~\ref{ass.prob} holds, $\mu_k = \mu \in \R{}_{>0}$ for all $k \in \N{}$, and there exists $\lambda \in \R{}_{>0}$ and infinite $\Scal_\lambda \subseteq \N{}$ such that for some $\{x_k\} \subset \Xcal \cap \R{n}_{\geq0}$ one finds $\nabla c(x_k)^T \nabla c(x_k) \succeq \lambda I$ for all $k \in \Scal_\lambda$.  Then, there exists $\kappa_{v,2} \in \R{}_{>0}$ such that
  \bequation\label{eq.v_2}
    \|c_k\|_2 - \|c_k + \nabla c(x_k)^T v_k\|_2 \geq \kappa_{v,2} \|v_k\|_2^2\ \ \text{for all}\ \ k \in \Scal_\lambda,
  \eequation
  where $v_k = u_k + \nabla c(x_k) w_k$ with $(u_k,w_k)$ being the unique optimal solution of \eqref{prob.v}.  Consequently, under the conditions of Lemma~\ref{lem.normal_stationarity} with $\Scal$ defined in that lemma, if $\Scal_\lambda$ defined as all sufficiently large indices in $\Scal$ satisfies the conditions above, then it follows that $\{v_k\}_{k \in \Scal} \to 0$ if and only if $\{\|c_k\|_2 - \|c_k + \nabla c(x_k)^T v_k\|_2\}_{k \in \Scal} \to 0$. 
\elemma
\bproof
  Consider arbitrary $k \in \Scal_\lambda$.  Under the stated conditions with $j_k := \nabla c(x_k)^T$, Lemma~\ref{lem.normal_generic} implies $\|c_k + j_k v_k\|_2 \leq \|c_k\|_2$.  Hence, by Assumption~\ref{ass.prob},
  \bequation\label{eq.vanilla}
    \baligned
      &\ \|c_k\|_2^2 - \|c_k + j_k v_k\|_2^2 = (\|c_k\|_2 + \|c_k + j_k v_k\|_2)(\|c_k\|_2 - \|c_k + j_k v_k\|_2) \\
      \leq&\ 2\|c_k\|_2(\|c_k\|_2 - \|c_k + j_k v_k\|_2) \leq 2\kappa_c(\|c_k\|_2 - \|c_k + j_k v_k\|_2).
    \ealigned
  \eequation
  If $v_k = 0$, then \eqref{eq.v_2} follows trivially.  Hence, we may proceed under the assumption that $v_k \neq 0$, which by $v_k = u_k + j_k^T w_k$ and the Fundamental Theorem of Linear Algebra means that $u_k \neq 0$ and/or $w_k \neq 0$.  If $w_k = 0$, then it follows by construction of \eqref{prob.v} that $u_k = 0$ as well.  Hence, we may conclude from $v_k \neq 0$ that, in fact, $w_k \neq 0$.  Since $(u_k,w_k)$ is the unique optimal solution of \eqref{prob.v}, it follows that $\alpha_k^* = 1$ is the optimal solution of the strongly convex quadratic optimization problem
  \bequation\label{prob.alpha01}
    \min_{\alpha\in[0,1]}\ \thalf\|c_k + \alpha j_k j_k^T w_k\|_2^2 + \thalf \mu_k \|\alpha u_k\|_2^2,
  \eequation
  which further implies (since an optimality condition of \eqref{prob.alpha01} is that the derivative of its objective function with respect to $\alpha$ is less than or equal to zero at $\alpha_k^* = 1$) that $-c_k^T j_k j_k^T w_k \geq \|j_k j_k^T w_k\|_2^2 + \mu_k\|u_k\|_2^2$.  Consequently, one finds
  \bequation\label{eq.intermediate}
    \baligned
      \|c_k\|_2^2 - \|c_k + j_k v_k\|_2^2
      &= \|c_k\|_2^2 - \|c_k + j_k j_k^T w_k\|_2^2 \\
      &= -2c_k^T j_k j_k^T w_k - \|j_k j_k^T w_k\|_2^2 \geq \|j_k j_k^T w_k\|_2^2 + 2\mu_k\|u_k\|_2^2.
    \ealigned
  \eequation
  With \eqref{eq.vanilla} and \eqref{eq.intermediate}, it follows from Assumption~\ref{ass.prob}, the conditions of the lemma, and since the Cauchy-Schwarz inequality implies $\|w_k\|_2 \geq \|j_k^T w_k\|_2/\|j_k^T\|_2$ that
  \begin{align*}
    \|c_k\|_2 - \|c_k + j_k v_k\|_2
    &\geq (2\kappa_c)^{-1} (\|c_k\|_2^2 - \|c_k + j_k v_k\|_2^2) \\
    &\geq (2\kappa_c)^{-1} (\|j_k j_k^T w_k\|_2^2 + 2\mu_k\|u_k\|_2^2) \\
    &\geq (2\kappa_c)^{-1} (\lambda^2 \|w_k\|_2^2 + 2\mu_k\|u_k\|_2^2) \\
    &\geq (2\kappa_c)^{-1} (\tfrac{\lambda^2}{\kappa_{\nabla c}^2} \|j_k^T w_k\|_2^2 + 2\mu_k\|u_k\|_2^2) \\
    &\geq (2\kappa_c)^{-1} \min\{\tfrac{\lambda^2}{\kappa_{\nabla c}^2},2\mu_k\} (\|j_k^T w_k\|_2^2 + \|u_k\|_2^2) \\
    &= (2\kappa_c)^{-1} \min\{\tfrac{\lambda^2}{\kappa_{\nabla c}^2},2\mu\} \|v_k\|_2^2 =: \kappa_{v,2} \|v_k\|_2^2,
  \end{align*}
  which gives \eqref{eq.v_2}, as desired.
\eproof

Next, we show that if the point defining subproblem~\eqref{prob.d} is not an infeasible stationary point for problem~\eqref{prob.opt}, then the subproblem with $g_k = \nabla f(x_k)$ yields a zero solution if and only if the point defining the subproblem is stationary for \eqref{prob.opt}.

\blemma\label{lem.stationarity_generic}
  Suppose that Assumption~\ref{ass.prob} holds and, with respect to $x \in \Xcal \cap \R{n}_{\geq0}$, one finds $c(x) = 0$.  Given $H \in \R{n \times n}$ with $H \succ 0$, consider $($recall \eqref{prob.d}$)$
  \bequation\label{prob.stationarity_generic}
    \min_{d\in\R{n}}\ \nabla f(x)^Td + \thalf d^THd\ \st\ c(x) + \nabla c(x)^T d = 0\ \ \text{and}\ \ x + d \geq 0.
  \eequation
  Then, one finds that the optimal solution of problem~\eqref{prob.stationarity_generic} is $d = 0$ if and only if $x$ is a KKT point $($i.e., first-order stationary point$)$ for problem~\eqref{prob.opt}.
\elemma
\bproof
  Suppose the conditions of the lemma hold and let $d$ be the optimal solution of \eqref{prob.stationarity_generic}.  Since $x \in \R{n}_{\geq0}$ and $c(x) = 0$, it follows that the zero vector is feasible for \eqref{prob.stationarity_generic}.  In addition, necessary and sufficient optimality conditions for subproblem~\eqref{prob.stationarity_generic} are that, corresponding to $d \in \R{n}$, there exist $y \in \R{m}$ and $z \in \R{n}$ such that
  \bequation\label{eq.KKT_fc_k}
    \nabla f(x) + Hd + \nabla c(x) y - z = 0,\ \ \nabla c(x)^T d = 0,\ \ \text{and}\ \ 0 \leq x + d \perp z \geq 0.
  \eequation
  If $d = 0$, then since $c(x)=0$ it follows that $(x,y,z)$ satisfies \eqref{eq.KKT}, as desired.  On the other hand, if~$x$ is a KKT point for~\eqref{prob.opt}, then there exist $y \in \R{m}$ and $z \in \R{n}$ such that $(x,y,z)$ satisfies \eqref{eq.KKT}, which in turn means that $d = 0$ along with $(y,z)$ satisfies \eqref{eq.KKT_fc_k}, and this solution is unique since the objective of \eqref{prob.stationarity_generic} is strongly convex.
\eproof

We conclude this subsection by showing that, under common assumptions and given $x_k \in \R{n}_{\geq0}$, the quantity $\|d_k\|_2^2$, where $d_k \in \R{n}$ solves subproblem~\eqref{prob.d} with $g_k = \nabla f(x_k)$, represents a stationarity measure with respect to \eqref{prob.opt}.   (The assumption in the lemma that $H_k = H$ for some $H \succ 0$ for all $k \in \N{}$ could be relaxed; see Remark~\ref{rem.mu_H} at the end of this subsection.  We consider this case for the sake of brevity.)

\blemma\label{lem.stationarity}
  Suppose that Assumption~\ref{ass.prob} holds and there exists infinite $\Scal \subseteq \N{}$ such that for some sequence $\{x_k\} \subset \Xcal \cap \R{n}_{\geq0}$ one finds $\{x_k\}_{k \in \Scal} \to x_*$ for some $x_* \in \Xcal \cap \R{n}_{\geq0}$ with $c(x_*) = 0$ and, with the notation in Lemma~\ref{lem.normal_stationarity}, one finds that
  \benumerate
    \item[(i)] $-\nabla f(x_*) = \nabla c(x_*) y - I_{\Acal(x_*)}^Tz_{\Acal(x_*)}$ for some $(y,z_{\Acal(x_*)}) \in \R{m} \times \R{|\Acal(x_*)|}_{>0}$ and
    \item[(ii)] the following matrix has full row rank:
    %\bequationNN
      $\bbmatrix \nabla c(x_*)^T \\ I_{\Acal(x_*)} \ebmatrix$.
    %\eequationNN
  \eenumerate
  Then, with $H_k = H$ for some $H \succ 0$ for all $k \in \N{}$, and with $d_k$ solving~\eqref{prob.d} with $g_k = \nabla f(x_k)$ for all $k \in \N{}$, $x_*$ satisfies~\eqref{eq.KKT} if and only if $\{\|d_k\|_2^2\}_{k \in \Scal} \to 0$.
\elemma
\bproof
  Letting $\Acal_* := \Acal(x_*)$ and considering the linear system of equations
  \bequationNN
    \bbmatrix H & \nabla c(x) & -I_{\Acal_*}^T \\ \nabla c(x)^T & 0 & 0 \\ I_{\Acal_*} & 0 & 0 \ebmatrix \bbmatrix d \\ y \\ z_{\Acal_*} \ebmatrix = \bbmatrix -\nabla f(x) \\ 0 \\ -x_{\Acal_*} \ebmatrix,
  \eequationNN
  the proof follows under the conditions of the lemma using the same line of deduction as the proof of Lemma~\ref{lem.normal_stationarity}, which we omit for the sake of brevity.
\eproof

\begin{remark}\label{rem.mu_H}
One might relax the condition in Lemma~\ref{lem.normal_stationarity} that $\mu = \mu_k$ for all $k \in \N{}$ and similarly relax the condition in Lemma~\ref{lem.stationarity} that $H_k = H \succ 0$ for all $k \in \N{}$, such as by requiring merely that $\{\mu_k\}_{k \in \Scal}$ and $\{H_k\}_{k \in \Scal}$ have bounded subsequences that converge to some $\mu \in \R{}_{>0}$ and $H \succ 0$, respectively.  In these cases, the ``if and only if'' statements would be replaced by an ``if'' statements, which in fact is all that is needed for our subsequent analysis and discussions.  Nevertheless, for brevity in the proofs, we provide the conditions that offer the stronger conclusions in these lemmas.
\end{remark}

%************
% Subsection
%************
\subsection{General Algorithm Behavior}\label{sec.general}

We now prove generally applicable results that hold for every run of Algorithm~\ref{alg.sqp}.  Our initial results in this section presume that iteration $k \in \N{}$ is reached and certain properties hold with respect to algorithmic quantities (e.g., $x_k \in \R{n}_{\geq0}$), although we ultimately prove in Lemma~\ref{lem.well_defined} that, in fact, these facts are guaranteed, i.e., they hold for any run for any generated $k \in \N{}$.  It is worthwhile to emphasize that the results in this section merely require that $g_k \in \R{n}$ for all $k \in \N{}$, which means, for example, that Assumption~\ref{ass.g} is not needed in this section.  All results that depend on the properties and effects of the stochastic gradient estimates are found in the subsequent subsection, i.e., Section~\ref{sec.guarantees}.

Our first lemma follows directly from Lemma~\ref{lem.normal_generic}, so it is stated without proof.

\blemma\label{lem.normal}
  Suppose that Assumption~\ref{ass.prob} holds.  Then, in any run of the algorithm such that iteration $k \in \N{}$ is reached and $x_k \in \R{n}_{\geq0}$, it holds that $v_k = 0$ if and only if $x_k$ satisfies \eqref{eq.KKT_con}, i.e., $x_k$ is either feasible or an infeasible stationary point, whereas $v_k \neq 0$ if and only if $\|c_k\|_2 > \|c_k + \nabla c(x_k)^Tv_k\|_2$.
\elemma

Our next result shows that, in any iteration in which the current iterate $x_k$ is in the nonnegative orthant and $\tau_{k-1} > 0$, the merit parameter is either kept at the same value or decreased, and, if it is decreased, then it is decreased below a constant fraction times its former value.  As in other SQP methods with such a feature, this ensures that if the merit parameter sequence does not vanish (i.e., its limiting value is nonzero), then it eventually remains at a constant positive value; see Lemma~\ref{lem.well_defined}.

\blemma\label{lem.merit_parameter_monotonic}
  Suppose that Assumption~\ref{ass.prob} holds. In any run of the algorithm such that line~\ref{step.g} of iteration $k \in \N{}$ is reached, $x_k \in \R{n}_{\geq0}$, and $\tau_{k-1} \in \R{}_{>0}$, it holds that $0 < \tau_k \leq \tau_{k-1}$, where if $\tau_k < \tau_{k-1}$, then $\tau_k \leq (1 - \epsilon_\tau) \tau_{k-1}$.
\elemma
\bproof
  Consider an arbitrary run in which line~\ref{step.g} of iteration $k \in \N{}$ is reached, $x_k \in \R{n}_{\geq0}$, and $\tau_{k-1} \in \R{}_{>0}$.  Let us show that $0 < \tau_k \leq \tau_{k-1}$, in which case the fact that $\tau_k < \tau_{k-1}$ implies $\tau_k \leq (1 - \epsilon_\tau) \tau_{k-1}$ follows from \eqref{eq.merit_parameter}.  Toward this end, let us next show that $\tautrial_k > 0$.  By the constraints of \eqref{prob.d}, \eqref{eq.tautrial}, and Lemma~\ref{lem.normal}, one finds that $\tautrial_k > 0$ whenever $\|c_k\|_2 - \|c_k + \nabla c(x_k)^Tv_k\|_2 > 0$.  Hence, to show that one always finds $\tautrial_k > 0$, all that remains is to consider the case when $\|c_k\|_2 - \|c_k + \nabla c(x_k)^Tv_k\|_2 = 0$.  In this case, it follows from Lemma~\ref{lem.normal} that $v_k = 0$, meaning that $d = 0$ is feasible for \eqref{prob.d}.  This, in turn, means that $g_k^Td_k + \thalf d_k^TH_kd_k \leq 0$, so by~\eqref{eq.tautrial} one finds that $\tautrial_k = \infty > 0$.  Since it has been shown that $\tautrial_k > 0$, the fact that $0 < \tau_k \leq \tau_{k-1}$ now follows directly from \eqref{eq.merit_parameter}, completing the proof.
\eproof

We now show that the model reduction offered by the computed search direction satisfies a lower bound with the properties stated in our algorithm development.

\blemma\label{lem.model_reduction}
  Suppose that Assumptions~\ref{ass.prob} and \ref{ass.H} hold.  In any run of the algorithm such that line~\ref{step.g} is reached in iteration $k \in \N{}$, $x_k \in \R{n}_{\geq0}$, and $\tau_k \in \R{}_{>0}$, one finds with $\zeta$ from Assumption~\ref{ass.H} that
  \bequation\label{eq.merit_reduction_lower}
    \Delta l(x_k,\tau_k,g_k,d_k) \geq \thalf \tau_k \zeta \|d_k\|_2^2 + \sigma(\|c_k\|_2 - \|c_k + \nabla c(x_k)^Td_k\|_2),
  \eequation
  and, if $d_k \neq 0$, then $\Delta l(x_k,\tau_k,g_k,d_k) > 0$.
\elemma
\bproof
  Consider an arbitrary run in which line~\ref{step.g} of iteration $k \in \N{}$ is reached, $x_k \in \R{n}_{\geq0}$, and $\tau_k \in \R{}_{>0}$.  By \eqref{eq.model_reduction} and Assumption~\ref{ass.H}, \eqref{eq.merit_reduction_lower} is implied by
  \bequation\label{eq.merit_reduction_lower_equiv}
    (1 - \sigma)(\|c_k\|_2 - \|c_k + \nabla c(x_k)^Td_k\|_2) \geq \tau_k (g_k^Td_k + \thalf d_k^TH_kd_k).
  \eequation
  If $g_k^Td_k + \thalf d_k^TH_kd_k \leq 0$, then \eqref{eq.merit_reduction_lower_equiv} holds due to Lemma~\ref{lem.normal} and the fact that \eqref{prob.d} ensures $\nabla c(x_k)^Tv_k = \nabla c(x_k)^Td_k$.  On the other hand, if $g_k^Td_k + \thalf d_k^TH_kd_k > 0$, then one finds by the update of the merit parameter, namely, \eqref{eq.tautrial} and \eqref{eq.merit_parameter}, that
  \bequationNN
    \tau_k \leq \tautrial_k = \tfrac{(1-\sigma)(\|c_k\|_2 - \|c_k + \nabla c(x_k)^Td_k\|_2) }{g_k^Td_k + \thalf d_k^TH_kd_k},
  \eequationNN
  from which \eqref{eq.merit_reduction_lower_equiv} follows again.  Finally, that $d_k \neq 0$ implies $\Delta l(x_k,\tau_k,g_k,d_k) > 0$ follows from \eqref{eq.merit_reduction_lower}, $\tau_k \in \R{}_{>0}$, and since $\zeta\in\R{}_{>0}$ in Assumption~\ref{ass.H}.
\eproof

Our next result is that, under the same conditions as our previous lemmas and under the assumption that $\xi_{k-1} \in \R{}_{>0}$, the ratio parameter is either kept at the same value or decreased, and, like the merit parameter, if it is decreased, then it is decreased at least below a constant fraction times its previous value.

\blemma\label{lem.xi}
  Suppose that Assumptions~\ref{ass.prob} and \ref{ass.H} hold.  In any run of the algorithm such that line~\ref{step.g} is reached in iteration $k \in \N{}$, $x_k \in \R{n}_{\geq0}$, $\tau_k \in \R{}_{>0}$, and $\xi_{k-1} \in \R{}_{>0}$, it holds that $0 < \xi_k \leq \xi_{k-1}$, where if $\xi_k < \xi_{k-1}$, then $\xi_k \leq (1 - \epsilon_\xi) \xi_{k-1}$.
\elemma
\bproof
  Consider an arbitrary run in which line~\ref{step.g} of iteration $k \in \N{}$ is reached, $x_k \in \R{n}_{\geq0}$, $\tau_k \in \R{}_{>0}$, and $\xi_{k-1} \in \R{}_{>0}$.  Let us show that $0 < \xi_k \leq \xi_{k-1}$, in which case the fact that $\xi_k < \xi_{k-1}$ implies $\xi_k \leq (1 - \epsilon_\xi) \xi_{k-1}$ follows from~\eqref{eq.xi_trial}.  Toward this end, observe that if $d_k = 0$, then the algorithm sets $\xi_k \gets \xi_{k-1} > 0$, which is consistent with the desired conclusion.  On the other hand, if $d_k \neq 0$, then by \eqref{eq.xi_trial}, $\tau_k \in \R{}_{>0}$, Lemma~\ref{lem.normal}, the fact that \eqref{prob.d} ensures $\nabla c(x_k)^Tv_k = \nabla c(x_k)^Td_k$, and Lemma~\ref{lem.model_reduction},
  \bequation\label{eq.xitrial_proof}
    \xitrial_k = \tfrac{\Delta l(x_k,\tau_k,g_k,d_k)}{\tau_k \|d_k\|_2^2} \geq \tfrac{\thalf \tau_k \zeta \|d_k\|_2^2}{\tau_k \|d_k\|_2^2} = \thalf \zeta > 0.
  \eequation
  Hence, by \eqref{eq.xi_trial}, the desired conclusion follows.
\eproof

Next, we prove bounds for the step size computed in the algorithm.

\blemma\label{lem.stoch_perf}
  Suppose that Assumptions~\ref{ass.prob} and \ref{ass.H} hold.  In any run of the algorithm such that line~\ref{step.g} is reached in iteration $k \in \N{}$, $x_k \in \R{n}_{\geq0}$, $\tau_k \in \R{}_{>0}$, and $\xi_k \in \R{}_{>0}$, it holds that $0 < \alpha_k^{\min} \leq \alpha_k^{\max} \leq \min\{1,\alpha_k^\varphi\}$, and, so, $x_{k+1} \in \R{n}_{\geq0}$.
\elemma
\bproof
  Consider an arbitrary run of the algorithm in which line~\ref{step.g} of iteration $k \in \N{}$ is reached, $x_k \in \R{n}_{\geq0}$, $\tau_k \in \R{}_{>0}$, and $\xi_k \in \R{}_{>0}$. Let us show that $0 < \alpha_k^{\min} \leq \alpha_k^{\max} \leq 1$, in which case the fact that $x_{k+1} \in \R{n}_{\geq0}$ follows from $x_k \in \R{n}_{\geq0}$, the fact that the constraints of \eqref{prob.d} ensure that $x_k + d_k \in \R{n}_{\geq0}$, and since the step size has $\alpha_k \in [\alpha_k^{\min}, \alpha_k^{\max}] \subset (0,1]$.  Toward this end, observe that if $d_k = 0$, then the algorithm yields $\alpha_k = \alpha_k^{\min} = \alpha_k^{\max} = \alpha_k^\varphi = 1$, so the conclusion follows trivially.  Hence, let us assume $d_k \neq 0$.  Observe that from \eqref{eq.beta_constraint}, the algorithm uses $\alpha_k^{\min}$ with
  \bequation\label{eq.alphamin_ub_1}
    0 < \alpha_k^{\min} = \tfrac{2 (1-\eta) \beta_k \xi_k \tau_k}{\tau_k L + \Gamma} \leq 1.
  \eequation
  Now observing \eqref{eq.alpha_max}, which shows $\alpha_k^{\max} \leq \min\{1,\alpha_k^\varphi\}$, one finds that all that remains is to prove that $\alpha_k^{\min} \leq \alpha_k^\varphi$.  For this purpose, let us introduce
  \bequationNN
    \alphasuff_k := \min\left\{ 1, \tfrac{2 (1-\eta) \beta_k \Delta l(x_k,\tau_k,g_k,d_k)}{(\tau_k L + \Gamma) \|d_k\|_2^2}\right\},
  \eequationNN
  where $\alphasuff_k \in (0,1]$ follows by $\beta_k \in (0,1]$, Lemma~\ref{lem.model_reduction}, and $d_k \neq 0$.  To show that $\alpha_k^{\min} \leq \alpha_k^\varphi$, our aim is to show that $\alpha_k^{\min} \leq \alphasuff_k \leq \alpha_k^\varphi$.  First, from \eqref{eq.xi_trial}, one finds
  \bequation\label{eq.alphamin_ub_2}
    \alpha_k^{\min} = \tfrac{2(1-\eta)\beta_k\xi_k\tau_k}{\tau_k L + \Gamma} \leq \tfrac{2(1-\eta)\beta_k\xitrial_k\tau_k}{\tau_k L + \Gamma} = \tfrac{2(1-\eta)\beta_k\Delta l(x_k,\tau_k,g_k,d_k)}{(\tau_kL + \Gamma)\|d_k\|_2^2}.
  \eequation
  Combining \eqref{eq.alphamin_ub_1} and \eqref{eq.alphamin_ub_2}, one finds that $\alpha_k^{\min} \leq \alphasuff_k$, as desired.  Now, toward proving that  $\alphasuff_k\leq\alpha_k^{\varphi}$, let us first show that $\varphi_k(\alphasuff_k)\leq 0$.  From the triangle inequality, the fact that $\alphasuff_k\in(0,1]$, and \eqref{eq.varphi}, it follows that
  \begin{align*}
    &\ \varphi_k(\alphasuff_k) \\
    =&\ (\eta - 1) \alphasuff_k \beta_k \Delta l(x_k,\tau_k,g_k,d_k) + \|c_k + \alphasuff_k \nabla c(x_k)^T d_k\|_2 - \|c_k\|_2 \\
    &\ + \alphasuff_k(\|c_k\|_2 - \|c_k + \nabla c(x_k)^Td_k\|_2) + \thalf (\tau_k L + \Gamma) (\alphasuff_k)^2 \|d_k\|_2^2 \\
    \leq&\ (\eta - 1) \alphasuff_k \beta_k \Delta l(x_k,\tau_k,g_k,d_k) + (1 - \alphasuff_k) \|c_k\|_2 + \alphasuff_k \|c_k + \nabla c(x_k)^Td_k\|_2 \\
    &\ - \|c_k\|_2 + \alphasuff_k(\|c_k\|_2 - \|c_k + \nabla c(x_k)^Td_k\|_2) + \thalf (\tau_k L + \Gamma) (\alphasuff_k)^2 \|d_k\|_2^2 \\
    =&\ (\eta - 1) \alphasuff_k \beta_k \Delta l(x_k,\tau_k,g_k,d_k) + \thalf (\tau_k L + \Gamma) (\alphasuff_k)^2 \|d_k\|_2^2 \\
    \leq&\ (\eta - 1) \alphasuff_k \beta_k \Delta l(x_k,\tau_k,g_k,d_k) + \thalf \alphasuff_k(\tau_k L + \Gamma) \|d_k\|_2^2 \(\tfrac{2(1-\eta) \beta_k \Delta l(x_k,\tau_k,g_k,d_k)}{(\tau_k L + \Gamma)\|d_k\|_2^2}\) \\
    =&\ 0.
  \end{align*}
  Therefore, by \eqref{eq.alpha_max}, it follows that $\alphasuff_k \leq \alpha_k^{\varphi}$.
\eproof

Our next lemma shows an upper bound on the change in the merit function.  In the lemma and throughout the rest of the paper, for any $k \in \N{}$ such that line~\ref{step.g} is reached we let $\dtrue_k \in \R{n}$ denote the solution of \eqref{prob.d} when $g_k$ is replaced by $\nabla f(x_k)$.

\blemma\label{lem.stoch_suff_decrease}
  Suppose that Assumptions~\ref{ass.prob} and \ref{ass.H} hold.  In any run of the algorithm such that line~\ref{step.g} is reached in iteration $k \in \N{}$, $x_k \in \R{}_{\geq0}$, $\tau_k \in \R{}_{>0}$, and $\alpha_k \in (0,\alpha_k^\varphi]$, it holds that
  \begin{multline*}
    \phi(x_k + \alpha_k d_k,\tau_k) - \phi(x_k,\tau_k) \leq -\alpha_k\Delta l(x_k,\tau_k,\nabla f(x_k),\dtrue_k) \\ + \alpha_k\tau_k\nabla f(x_k)^T(d_k - \dtrue_k) + (1-\eta)\alpha_k\beta_k\Delta l(x_k,\tau_k,g_k,d_k).
  \end{multline*}
\elemma
\bproof
  Consider an arbitrary run of the algorithm in which line~\ref{step.g} of iteration $k \in \N{}$ is reached, $x_k \in \R{n}_{\geq0}$, $\tau_k \in \R{}_{>0}$, and $\alpha_k \in (0,\alpha_k^\varphi]$.  By Assumption~\ref{ass.prob} (which led to \eqref{eq.Lipschitz}), \eqref{prob.d} (which implies $c_k + \nabla c(x_k)^Td_k = c_k + \nabla c(x_k)^T\dtrue_k$), \eqref{eq.model_reduction}, \eqref{eq.varphi}, and the fact that $0 < \alpha_k \leq \alpha_k^\varphi$ (which means $\varphi_k(\alpha_k) \leq 0$), it follows that
  \begin{align*}
    &\ \phi(x_k + \alpha_k d_k,\tau_k) - \phi(x_k,\tau_k) \\
    =&\ \tau_k(f(x_k + \alpha_k d_k) - f_k) + \|c(x_k + \alpha_k d_k)\|_2 - \|c_k\|_2 \\
    \leq&\ \alpha_k\tau_k\nabla f(x_k)^Td_k + \|c_k + \alpha_k\nabla c(x_k)^Td_k\|_2 - \|c_k\|_2 + \thalf(\tau_k L + \Gamma)\alpha_k^2\|d_k\|_2^2 \\
    =&\ -\alpha_k\Delta l(x_k,\tau_k,\nabla f(x_k),\dtrue_k) + \alpha_k\tau_k\nabla f(x_k)^T(d_k - \dtrue_k) + \|c_k + \alpha_k\nabla c(x_k)^Td_k\|_2 \\
    &\ - \|c_k\|_2 + \alpha_k(\|c_k\|_2 - \|c_k + \nabla c(x_k)^Td_k\|_2) + \thalf(\tau_k L + \Gamma)\alpha_k^2\|d_k\|_2^2 \\
    \leq&\ -\alpha_k\Delta l(x_k,\tau_k,\nabla f(x_k),\dtrue_k) + \alpha_k\tau_k\nabla f(x_k)^T(d_k - \dtrue_k) \\
    &\ + (1-\eta)\alpha_k\beta_k\Delta l(x_k,\tau_k,g_k,d_k),
  \end{align*}
  which shows the desired conclusion.
\eproof

We now show that each search direction---and, similarly, the search direction that would be computed if the true gradient of the objective function were used in place of the stochastic gradient estimate---can be viewed as a projection of the unconstrained minimizer of the objective of~\eqref{prob.d} onto the feasible region of~\eqref{prob.d}.

\blemma\label{lem.stationarity_projection}
  Suppose that Assumptions~\ref{ass.prob} and \ref{ass.H} hold.  In any run of the algorithm such that line~\ref{step.g} is reached in iteration $k \in \N{}$, $x_k \in \R{}_{\geq0}$, and with
  \bequationNN
    \Dcal_k := \{d \in \R{n} : \nabla c(x_k)^T(d - v_k ) = 0 , \ x_k + d \geq 0\}\ \text{and}\ 
    \proj_k(\dbar) := \argmin_{d \in \Dcal_k} \|d - \dbar\|_{H_k}^2,
  \eequationNN
  it holds that $d_k = \proj_k(-H_k^{-1}g_k)$ and $\dtrue_k = \proj_k(-H_k^{-1} \nabla f(x_k))$.
\elemma
\bproof
  Consider an arbitrary run of the algorithm in which line~\ref{step.g} of iteration $k \in \N{}$ is reached and $x_k \in \R{}_{\geq0}$.  The desired conclusion follows from the facts that $\Dcal_k$ is convex and, under Assumption~\ref{ass.H}, $H_k$ is SPD; in particular, one finds that
  \begin{align*}
    d_k &= \argmin_{d\in\Dcal_k}\ g_k^Td + \thalf d^TH_kd = \argmin_{d\in\Dcal_k}\ \thalf\|d + H_k^{-1}g_k\|_{H_k}^2 = \proj_k(-H_k^{-1}g_k),
  \end{align*}
  and similarly with respect to $\dtrue_k$ with $g_k$ replaced by $\nabla f(x_k)$.
\eproof

We are now prepared to prove the following lemma, which shows that the algorithm is well defined and either terminates finitely with an infeasible stationary point or generates an infinite sequence of iterates with certain critical properties of the simultaneously generated algorithmic sequences. The lemma also reveals that the monotonically nonincreasing merit parameter sequence either vanishes or ultimately remains constant, and it reveals that the monotonically nonincreasing ratio parameter sequence ultimately remains constant at a value that is greater than or equal to a positive real number that is defined uniformly across all runs of the algorithm.

\blemma\label{lem.well_defined}
  Suppose that Assumptions~\ref{ass.prob} and \ref{ass.H} hold.  In any run, either the algorithm terminates finitely with an infeasible stationary point or it performs an infinite number of iterations such that, for all $k \in \N{}$, it holds that
  \benumerate
    \item[(a)] $x_k \in \R{n}_{\geq0}$,
    \item[(b)] $v_k = 0$ if and only if $x_k$ satisfies \eqref{eq.KKT_con},
    \item[(c)] $v_k \neq 0$ if and only if $\|c_k\|_2 > \|c_k + \nabla c(x_k)^Tv_k\|_2$,
    \item[(d)] $0 < \tau_k \leq \tau_{k-1} < \infty$,
    \item[(e)] $\tau_k < \tau_{k-1}$ if and only if $\tau_k \leq (1 - \epsilon_\tau) \tau_{k-1}$,
    \item[(f)] \eqref{eq.merit_reduction_lower} holds,
    \item[(g)] $d_k \neq 0$ if and only if $\Delta l(x_k,\tau_k,g_k,d_k) > 0$,
    \item[(h)] $0 < \xi_k \leq \xi_{k-1} < \infty$,
    \item[(i)] $\xi_k < \xi_{k-1}$ if and only if $\xi_k \leq (1 - \epsilon_\xi) \xi_{k-1}$, and
    \item[(j)] $0 < \alpha_k^{\min} \leq \alpha_k^{\max} \leq \min\{1, \alpha_k^\varphi\}$.
  \eenumerate
  In addition, in any run that does not terminate finitely, it holds that
  \benumerate
    \item[(k)] either $\{\tau_k\} \searrow 0$ or there exists $k_\tau \in \N{}$ and $\tau_{\min} \in \R{}_{>0}$ such that $\tau_k = \tau_{\min}$ for all $k \in \N{}$ with $k \geq k_\tau$, and
    \item[(l)] there exist $k_\xi \in \N{}$ and $\xi_{\min} \in \R{}_{>0}$ with $\xi_{\min} \geq \thalf \zeta (1 - \epsilon_\xi)$ such that $\xi_k = \xi_{\min}$ for all $k \in \N{}$ with $k \geq k_\xi$.
  \eenumerate
\elemma
\bproof
  Given the initialization of the algorithm, statements $(a)$--$(j)$ follow by induction from Lemmas~\ref{lem.normal}--\ref{lem.stoch_perf}.  Statement $(k)$ follows from statements $(d)$ and~$(e)$.  Finally, to prove statement $(l)$, consider arbitrary $k \in \N{}$ in a run that does not terminate finitely and note that if $d_k = 0$, then $\xitrial_k \gets \infty$, and if $d_k \neq 0$, then $\xitrial_k$ satisfies \eqref{eq.xitrial_proof}, meaning that $\xitrial_k \geq \thalf \zeta$.  Consequently, by~\eqref{eq.xi_trial}, $\xi_k < \xi_{k-1}$ only if $\xi_{k-1} > \thalf \zeta$.  This, along with statements $(h)$ and $(i)$, leads to the conclusion.
\eproof

%************
% Subsection
%************
\subsection{Convergence Guarantees}\label{sec.guarantees}

We now turn to prove convergence results under Assumption~\ref{ass.good_merit_parameter} below.  Recalling the role of $\thalf \zeta (1 - \epsilon_\xi) \in \R{}_{>0}$ in Lemma~\ref{lem.well_defined}$(l)$, the assumption focuses on the following event for some $(k_{\min},\tau_{\min},f_{\sup}) \in \N{} \times \R{}_{>0} \times \R{}$, where for all generated $k \in \N{}$ we denote $\tautruetrial_k$ as the value of $\tautrial_k$ that would be computed in iteration $k$ if \eqref{prob.d} were solved with $\nabla f(x_k)$ in place of $g_k$:
\vspace{-15pt}

\bequation\label{def.E}
\baligned
  & \Ecal(k_{\min},\tau_{\min},f_{\sup}) \\
  := \{ &\text{An infinite number of iterations are performed, $f(x_{k_{\min}}) \leq f_{\sup}$, and} \\
  & \text{there exist $k' \in \N{}$ with $k' \leq k_{\min}$, $\tau' \in \R{}_{>0}$ with $\tau' \geq \tau_{\min}$,} \\
  & \text{and $\xi' \in \R{}_{>0}$ with $\xi' \geq \thalf \zeta (1 - \epsilon_\xi)$ such that} \\
  & \text{$\tau_k = \tau' \leq \tautruetrial_k$ and $\xi_k = \xi'$ for all $k \in \N{}$ with $k \geq k'$}\}.
\ealigned
\eequation

The following assumption is made in this subsection.  We present a discussion and supporting theoretical results about this assumption in Section~\ref{sec.non-vanishing}.

\bassumption\label{ass.good_merit_parameter}
  For some $(k_{\min},\tau_{\min},f_{\sup}) \in \N{} \times \R{}_{>0} \times \R{}$, the event $\Ecal := \Ecal(k_{\min}, \tau_{\min},f_{\sup})$ occurs and, conditioned on the occurrence of~$\Ecal$, Assumption~\ref{ass.prob} holds $($with the same constants as previously presented in \eqref{eq.bounds} and \eqref{eq.Lipschitz}$)$.
\eassumption

It is not a shortcoming of our analysis that Assumption~\ref{ass.good_merit_parameter}, through the definition of $\Ecal$, assumes that $(i)$ an infinite number of iterations are performed, $(ii)$ the objective value is bounded above in iteration $k_{\min}$, and $(iii)$ $\{\xi_k\}$ ultimately becomes a constant sequence with value at least $\thalf \zeta (1 - \epsilon_\xi) \in \R{}_{>0}$.  After all: $(i)$ Lemma~\ref{lem.well_defined} shows that the only alternative to an infinite number of iterations being performed is that the algorithm terminates finitely with an infeasible stationary point, in which case there is nothing else to prove; $(ii)$ $f_{\sup} \in \R{}$ can be arbitrarily large and knowledge of it is not required by the algorithm, so assuming that it exists is a very loose requirement; and $(iii)$ Lemma~\ref{lem.well_defined}$(l)$ shows that, in any run that does not terminate finitely, $\{\xi_k\}$ is monotonically nonincreasing and bounded below by $\thalf \zeta (1 - \epsilon_\xi) \in \R{}_{>0}$, which is a constant, i.e., it is not run-dependent.  Overall, the only important restriction of our analysis in this section is the fact that $\Ecal$ includes the requirement that $\{\tau_k\}$ ultimately becomes constant at a value at least $\tau_{\min}$ that is sufficiently small relative to $\{\tautruetrial_k\}$.  This restriction is the subject of Section~\ref{sec.non-vanishing}.

For the remainder of this subsection, we consider the stochastic process corresponding to the statement of Algorithm~\ref{alg.sqp}.  Specifically, the sequence
\bequationNN
  \{(x_k, v_k, g_k, d_k, \dtrue_k, \tautrial_k, \tautruetrial_k, \tau_k, \xitrial_k, \xi_k, \alpha_k^{\min}, \alpha_k^\varphi, \alpha_k^{\max}, \alpha_k)\}
\eequationNN
generated in any run can be viewed as a realization of the stochastic process
\bequationNN
  \{(X_k, V_k, G_k, D_k, \Dtrue_k, \Tautrial_k, \Tautruetrial_k, \Tau_k, \Xitrial_k, \Xi_k, \Alpha_k^{\min}, \Alpha_k^\varphi, \Alpha_k^{\max}, \Alpha_k)\}.
\eequationNN
Let $\Gcal_1$ denote the $\sigma$-algebra defined by the initial conditions of the algorithm and, for all $k \in \N{}$ with $k \geq 2$, let $\Gcal_k$ denote the $\sigma$-algebra generated by the initial conditions and the random variables $\{G_1,\dots,G_{k-1}\}$.  Then, with respect to the event $\Ecal$ in Assumption~\ref{ass.good_merit_parameter}, denote the trace $\sigma$-algebra of $\Ecal$ on $\Gcal_k$ as $\Fcal_k := \Gcal_k \cap \Ecal$ for all $k \in \N{}$.  It follows that $\{\Fcal_k\}$ is a filtration, and we proceed in our analysis under Assumptions~\ref{ass.g}, \ref{ass.H}, and \ref{ass.good_merit_parameter} (which subsumes Assumption~\ref{ass.prob}) with the definitions 
\bequationNN
  \P_k[\cdot] := \P_\omega[\cdot | \Fcal_k]\ \ \text{and}\ \ \E_k[\cdot] := \E_\omega[\cdot | \Fcal_k]
\eequationNN
(where $P_\omega$ denotes probability taken with respect to the distribution of $\omega$).  We also define, with respect to $\Ecal$, the random variables $K' \leq k_{\min}$, $\Tau' \geq \tau_{\min}$, and $\Xi' \geq \thalf \zeta (1 - \epsilon_\xi)$, which for a given run of the algorithm have the realized values $k'$, $\tau'$, and~$\xi'$, respectively, defined in \eqref{def.E}.  Conditioned on $\Ecal$, one has in any run that
\bequation\label{eq.tau_xi_random}
  \tau_{\min} \leq \Tau' \leq \tau_0\ \ \text{and}\ \ \thalf \zeta (1 - \epsilon_\xi) \leq \Xi' \leq \xi_0,
\eequation
and one has that $\Tau'$ and $\Xi'$ are $\Fcal_k$-measurable for $k = k_{\min} \geq K'$.

Our next lemma shows upper bounds on the norm of the difference between the computed search direction and the search direction that would be computed with the true gradient of the objective.  (The conclusion of this lemma and the following one would hold even without assuming that the event $\Ecal$ occurs, but in each result we condition on $\Fcal_k := \Gcal_k \cap \Ecal$ for use in our ultimate results under $\Ecal$.)

\blemma\label{lem.d_diff_bound}
  Suppose that Assumptions~\ref{ass.g}, \ref{ass.H}, and \ref{ass.good_merit_parameter} hold.  For all $k \in \N{}$,
  \begin{align*}
    \|D_k - \Dtrue_k\|_2 &\leq \zeta^{-1} \|G_k - \nabla f(X_k)\|_2 \\ \text{and}\ \ 
    \E_k[\|D_k - \Dtrue_k\|_2] &\leq \zeta^{-1} \E_k[\|G_k - \nabla f(X_k)\|_2] \leq \zeta^{-1} \sqrt{\rho_k}.
  \end{align*}
\elemma
\bproof
  Consider arbitrary $k \in \N{}$ under the stated conditions.  Lemma~\ref{lem.stationarity_projection} and the obtuse angle lemma for projections~\cite[Proposition~1.1.9]{Bert09} imply
  \begin{align*}
    (D_k - \Dtrue_k)^T H_k (-H_k^{-1}\nabla f(X_k) - \Dtrue_k) &\leq 0 \\ \text{and}\ \ 
    (\Dtrue_k - D_k)^T H_k (-H_k^{-1}G_k - D_k) &\leq 0.
  \end{align*}
  Summing these inequalities yields
  \begin{align*}
    0 &\geq (D_k - \Dtrue_k)^TH_k(-H_k^{-1}\nabla f(X_k) - \Dtrue_k) + (\Dtrue_k - D_k)^TH_k(-H_k^{-1}G_k - D_k) \\
    %&= (D_k - \Dtrue_k)^TH_k(-H_k^{-1}\nabla f(X_k) - \Dtrue_k + H_k^{-1}G_k + D_k) \\
    &= \|D_k - \Dtrue_k\|_{H_k}^2 - (D_k - \Dtrue_k)^T(\nabla f(X_k) - G_k).
  \end{align*}
  Hence, by the Cauchy–Schwarz inequality, it follows that
  \bequationNN
    \|D_k - \Dtrue_k\|_{H_k}^2 \leq (D_k - \Dtrue_k)^T(\nabla f(X_k) - G_k) \leq \|D_k - \Dtrue_k\|_2\|\nabla f(X_k) - G_k\|_2,
  \eequationNN
  which shows under Assumption~\ref{ass.H} that $\|D_k - \Dtrue_k\|_2 \leq \zeta^{-1} \|G_k - \nabla f(X_k)\|_2$, as desired.  Then, from this inequality, Assumption~\ref{ass.g}, and Jensen's inequality, one has
  \bequationNN
    \E_k[\|G_k - \nabla f(X_k)\|_2] \leq \sqrt{\E_k[\|G_k - \nabla f(X_k)\|_2^2]} \leq \sqrt{\rho_k},
  \eequationNN
  from which the remainder of the conclusion follows.
\eproof

We now show an upper bound on the expected difference between inner products involving the true and stochastic gradients and the true and stochastic directions.

\blemma\label{lem.gd_diff}
  Suppose that Assumptions~\ref{ass.g}, \ref{ass.H}, and \ref{ass.good_merit_parameter} hold.  For all $k \geq k_{\min}$,
  \begin{align*}
    |\E_k[G_k^TD_k - \nabla f(X_k)^T\Dtrue_k]| &\leq \zeta^{-1} (\rho_k + \kappa_{\nabla f}\sqrt{\rho_k}) \\ \text{and}\ \ 
    \E_k[\Delta l(X_k,\Tau_k,G_k,D_k)] - \Delta l(X_k,\Tau',\nabla f(X_k),\Dtrue_k) &\leq \Tau' \zeta^{-1} (\rho_k + \kappa_{\nabla f} \sqrt{\rho_k}).
  \end{align*}
\elemma
\bproof
  Consider arbitrary $k \geq k_{\min}$ under the stated conditions.  From the triangle and Cauchy–Schwarz inequalities and Lemma~\ref{lem.d_diff_bound}, it holds that
  \begin{align*}
    &\ |\E_k[G_k^TD_k - \nabla f(X_k)^T\Dtrue_k]| \\
    =&\ |\E_k[(G_k - \nabla f(X_k))^T\Dtrue_k + (G_k-\nabla f(X_k))^T(D_k - \Dtrue_k) \\
    &\ + \nabla f(X_k)^T(D_k - \Dtrue_k)]| \\
    =&\ |\E_k[(G_k-\nabla f(X_k))^T(D_k - \Dtrue_k)] + \E_k[\nabla f(X_k)^T(D_k - \Dtrue_k)]| \\
    \leq&\ \E_k[\|G_k - \nabla f(X_k)\|_2 \|D_k - \Dtrue_k\|_2] + \|\nabla f(X_k)\|_2 \E_k[\|D_k - \Dtrue_k\|_2] \\
    \leq&\ \zeta^{-1} \E_k[\|G_k - \nabla f(X_k)\|_2^2] + \zeta^{-1} \kappa_{\nabla f} \E_k[\|G_k - \nabla f(X_k)\|_2] \leq \zeta^{-1} \rho_k + \zeta^{-1} \kappa_{\nabla f} \sqrt{\rho_k},
  \end{align*}
  which gives the first result.  Then, for $k \geq k_{\min}$, \eqref{eq.model_reduction} and the equation above give
  \begin{align*}
    &\ \E_k[\Delta l(X_k,\Tau_k,G_k,D_k)] - \Delta l(X_k,\Tau',\nabla f(X_k),\Dtrue_k) \\
    =&\ \Tau' \E_k[\nabla f(X_k)^T\Dtrue_k - G_k^TD_k] \leq \Tau' \zeta^{-1} (\rho_k + \kappa_{\nabla f} \sqrt{\rho_k}),
  \end{align*}
  which completes the proof.
\eproof

Our next lemma shows a lower bound on the true model reduction.  In the lemma and our subsequent results, we define $J_k := \nabla c(X_k)^T$ for the sake of brevity.

\blemma\label{lem.daniel}
  Suppose that Assumptions~\ref{ass.g}, \ref{ass.H}, and \ref{ass.good_merit_parameter} hold.  For all $k \geq k_{\min}$,
  \bequationNN
    \Delta l(X_k,\Tau_k,\nabla f(X_k),\Dtrue_k) \geq \thalf \Tau' \zeta \|\Dtrue_k\|_2^2 + \sigma(\|c(X_k)\|_2 - \|c(X_k) + J_k\Dtrue_k\|_2) \geq 0.
  \eequationNN
\elemma
\bproof
  Consider arbitrary $k \geq k_{\min}$ under the stated conditions.  By \eqref{eq.model_reduction}, the fact that $\Tau_k = \Tau'$, and Assumption~\ref{ass.H}, the first desired conclusion is implied by
  \bequationNN
    (1 - \sigma)(\|c(X_k)\|_2 - \|c(X_k) + J_k\Dtrue_k\|_2) \geq \Tau' (\nabla f(X_k)^T\Dtrue_k + \thalf (\Dtrue_k)^TH_k \Dtrue_k).
  \eequationNN
  If $\nabla f(X_k)^T\Dtrue_k + \thalf (\Dtrue_k)^TH_k\Dtrue_k \leq 0$, then the above holds due to Lemma~\ref{lem.well_defined} and the fact that $J_k\Dtrue_k = J_kV_k$; else, $\nabla f(X_k)^T\Dtrue_k + \thalf (\Dtrue_k)^TH_k \Dtrue_k > 0$, in which case one finds from the conditions of the lemma, \eqref{eq.tautrial}, and \eqref{eq.merit_parameter} that
  \bequationNN
    \Tau_k = \Tau' \leq \Tautruetrial_k = \tfrac{(1 - \sigma)(\|c(X_k)\|_2 - \|c(X_k) + J_k \Dtrue_k\|_2)}{\nabla f(X_k)^T\Dtrue_k + \thalf (\Dtrue_k)^TH_k \Dtrue_k},
  \eequationNN
  from which the displayed inequality above follows again.  Finally, the remaining desired conclusion follows from the first conclusion, Lemma~\ref{lem.well_defined}, and $J_k\Dtrue_k = J_kV_k$.
\eproof

Next, we prove a critical upper bound on the expected value of the second term on the right-hand side of the upper bound proved in Lemma~\ref{lem.stoch_suff_decrease}.

\blemma\label{lem.stoch_key_diff}
  Suppose that Assumptions~\ref{ass.g}, \ref{ass.H}, and \ref{ass.good_merit_parameter} hold.  For all $k \geq k_{\min}$,
  \bequationNN
    \E_k[\Alpha_k\Tau_k\nabla f(X_k)^T(D_k - \Dtrue_k)] \leq \left(\tfrac{2(1-\eta) \Xi' \Tau'}{\Tau'L+\Gamma} + \theta \right) \beta_k \Tau' \kappa_{\nabla f} \zeta^{-1} \sqrt{\rho_k}.
  \eequationNN
\elemma
\bproof
  For arbitrary $k \geq k_{\min}$ under the conditions, \eqref{eq.beta_constraint} and \eqref{eq.alpha_max} yield
  \bequation\label{eq.Alpha'}
    \Alpha_k^{\min} = \beta_k \Alpha'\ \ \text{and}\ \ \Alpha_k^{\max} \leq \Alpha_k^{\min} + \theta\beta_k,\ \ \text{where}\ \ \Alpha' = \tfrac{2(1-\eta) \Xi' \Tau'}{\Tau' L+\Gamma}.
  \eequation
  Letting $\Pcal_k$ denote the event that $\nabla f(X_k)^T(D_k - \Dtrue_k) \geq 0$ and letting $\Pcal_k^c$ denote the event that $\nabla f(X_k)^T(D_k - \Dtrue_k) < 0$, the law of total expectation and the fact that $\Tau'$ and $\Xi'$ are $\Fcal_k$-measurable for $k \geq k_{\min}$ shows that
  \begin{align*}
    &\ \E_k[\Alpha_k \Tau_k \nabla f(X_k)^T(D_k - \Dtrue_k)] \\
    =&\ \P_k[\Pcal_k] \cdot \E_k[\Alpha_k \Tau' \nabla f(X_k)^T(D_k - \Dtrue_k) | \Pcal_k] \\
    &\ + \P_k[\Pcal_k^c] \cdot \E_k[\Alpha_k \Tau' \nabla f(X_k)^T(D_k - \Dtrue_k) | \Pcal_k^c] \\
    \leq&\ (\Alpha_k^{\min} + \theta \beta_k) \Tau' \P_k[\Pcal_k] \cdot \E_k[\nabla f(X_k)^T(D_k - \Dtrue_k) | \Pcal_k] \\
    &\ + \Alpha_k^{\min} \Tau' \P_k[\Pcal_k^c] \cdot \E_k[\nabla f(X_k)^T(D_k - \Dtrue_k) | \Pcal_k^c] \\
    =&\ \Alpha_k^{\min} \Tau' \E_k[\nabla f(X_k)^T(D_k - \Dtrue_k)] \\
    &\ + \theta \beta_k \Tau' \P_k[\Pcal_k] \cdot \E_k[\nabla f(X_k)^T(D_k - \Dtrue_k) | \Pcal_k].
  \end{align*}
  The Cauchy-Schwarz inequality and law of total expectation show that
  \begin{align*}
    &\ \P_k[\Pcal_k] \cdot \E_k[\nabla f(X_k)^T(D_k - \Dtrue_k) | \Pcal_k] \\
    \leq&\ \P_k[\Pcal_k] \cdot \E_k[\|\nabla f(X_k)\|_2 \|D_k - \Dtrue_k\|_2 | \Pcal_k] \\
    =&\ \E_k[\|\nabla f(X_k)\|_2 \|D_k - \Dtrue_k\|_2] - \P_k[\Pcal_k^c] \cdot \E_k[\|\nabla f(X_k)\|_2 \|D_k - \Dtrue_k\|_2 | \Pcal_k^c] \\
    \leq&\ \E_k[\|\nabla f(X_k)\|_2 \|D_k - \Dtrue_k\|_2],
  \end{align*}
  so from above, the Cauchy-Schwarz inequality, Assumption~\ref{ass.good_merit_parameter}, and Lemma~\ref{lem.d_diff_bound},
  \begin{align*}
    \E_k[\Alpha_k \Tau_k \nabla f(X_k)^T(D_k - \Dtrue_k)] &\leq (\Alpha_k^{\min} + \theta \beta_k) \Tau' \|\nabla f(X_k)\|_2 \E_k[\|D_k - \Dtrue_k\|_2] \\
    &\leq \left(\tfrac{2(1-\eta)\Xi'\Tau'}{\Tau'L+\Gamma} + \theta \right) \beta_k \Tau' \kappa_{\nabla f} \zeta^{-1} \sqrt{\rho_k},
  \end{align*}
  which gives the desired conclusion.
\eproof

We now present, as a lemma, results pertaining to the asymptotic behavior of the model reductions generated by the algorithm.  In the subsequent theorem after the lemma, these results will be translated in terms of quantities that, as seen in Section~\ref{sec.subproblems}, can be connected to stationarity measures related to problem~\eqref{prob.opt}.  We remark that the conditions of the lemma can be satisfied in a run-dependent manner if, every time the merit or ratio parameter is decreased, say in iteration $\khat \in \N{}$, the sequence $\{\beta_k\}$ is ``restarted'' such that with $\alpha' = 2(1-\eta)\xi_{\khat}\tau_{\khat}/(\tau_{\khat} L + \Gamma)$ and some (run-independent) $\psi \in (0,1]$ one chooses $\beta_k = \beta = \psi \frac{\alpha'}{2(1-\eta)(\alpha' + \theta)}$ for part (a) of the lemma and $\beta_k = \frac{1}{k - \khat + 1} \psi \frac{\alpha'}{2(1-\eta)(\alpha' + \theta)}$ for part (b); such a scheme was described in \cite{BeraCurtRobiZhou21} as well.  Notice that in this situation, $\beta$ and $\{\beta_k\}_{k \geq \khat}$ in parts (a) and (b), respectively, are random variables, \emph{but, importantly, they are $\Fcal_k$-measurable for $k \geq k_{\min}$}.  Alternatively, one could choose $\{\beta_k\}$ using the same formulas, but with $\xi_{\min}$ and $\tau_{\min}$ in place of $\xi_k$ and $\tau_k$, respectively, in the formula for $\alpha'$, in which case the choices are run-independent.  The downside of relying on this latter situation is that it requires knowledge of $\xi_{\min}$ and $\tau_{\min}$, which would not typically be known \emph{a priori}.  Hence, we analyze the former scheme, but use run-dependent bounds that, under $\Ecal$, are defined with respect to $\xi_{\min}$ and $\tau_{\min}$ (even though these values are unknown).

We also remark that for case (a) in the following lemma, the sequence $\{\rho_k\}$, which bounds the expected squared error in the stochastic gradient estimates, can be a constant sequence.  However, for case (b), the relationship between $\{\rho_k\}$ and $\{\beta_k\}$ means that the expected squared error in the gradient estimates must vanish as $k \to \infty$.  This requirement, which is stronger than the requirement for equality-constraints-only case in \cite{BeraCurtRobiZhou21}, is needed to overcome the fact that in the presence of bound constraints the search directions can be biased estimates of their true counterparts.

\blemma\label{lem.stoch_main}
  Under Assumptions~\ref{ass.g},~\ref{ass.H}, and~\ref{ass.good_merit_parameter}, suppose that $\{\rho_k\}$ is chosen such that there exists $\iota \in \R{}_{>0}$ with $\rho_k \leq \iota \beta_k^2$ for all $k \in \N{}$ with $k \geq k_{\min}$, and define
  \begin{align*}
    \alpha_{\min}' &= \tfrac{2(1-\eta)\xi_{\min}\tau_{\min}}{\tau_{\min}L + \Gamma},\ \ \alpha_{\max}' = \tfrac{2(1-\eta)\xi_0\tau_0}{\tau_0L + \Gamma}, \\ \text{and}\ \ \rho_{\max}' &= (\alpha_{\max}' + \theta) \tau_0 \zeta^{-1} (\kappa_{\nabla f} \sqrt{\iota} + (1-\eta) (\iota + \kappa_{\nabla f} \sqrt{\iota})).
  \end{align*}
  Then, with $\Alpha'$ defined in \eqref{eq.Alpha'} and $\E[\cdot | \Ecal]$ denoting total expectation over all realizations of the algorithm conditioned on event $\Ecal$, the following statements hold true.
  \benumerate
    \item[(a)] if $\beta_k = \beta = \psi \frac{\Alpha'}{2(1-\eta)(\Alpha' + \theta)}$ for some $\psi \in (0,1]$ for all $k \geq k_{\min}$, then
  \eenumerate
  \bequationNN
    \limsup_{k\to\infty}\ \E\left[\frac{1}{k} \sum_{j=k_{\min}}^{k_{\min} + k - 1} \Delta l(X_j,\Tau',\nabla f(X_j),\Dtrue_j) \Bigg| \Ecal \right] \leq \tfrac{\psi (\alpha_{\max}')^2 (\alpha_{\min}' + \theta) \rho_{\max}'}{2 (1-\eta) (1 - \frac\psi2) (\alpha_{\min}')^2 (\alpha_{\max}' + \theta)^2};
  \eequationNN
  \benumerate
  \item[(b)] if $\sum_{k=k_{\min}}^{\infty}\beta_k = \infty$, $\sum_{k=k_{\min}}^{\infty}\beta_k^2 < \infty$, and $\beta_k \leq \psi\frac{\Alpha'}{2(1-\eta)(\Alpha' + \theta)}$ for some $\psi \in (0,1]$ for all $k \geq k_{\min}$, it holds that
  \eenumerate
  \bequationNN
    \E\left[ \frac{1}{\sum_{j = k_{\min}}^{k_{\min} + k - 1} \beta_j} \sum_{j=k_{\min}}^{k_{\min} + k - 1} \beta_j \Delta l(X_j,\Tau',\nabla f(X_j),\Dtrue_j) \Bigg| \Ecal \right] \xrightarrow{k \to \infty} 0.
  \eequationNN
\elemma
\bproof
  For arbitrary $k \geq k_{\min}$ under the conditions, it follows from Lemma~\ref{lem.stoch_suff_decrease}, Lemma~\ref{lem.daniel} (which shows $\Delta l(X_k, \Tau_k, \nabla f(X_k), \Dtrue_k) \geq 0$), \eqref{eq.Alpha'}, the fact that $\Alpha_k \geq \Alpha_k^{\min} = \Alpha' \beta_k$, Lemma~\ref{lem.stoch_key_diff}, the fact that $\Alpha_k \leq \Alpha_k^{\max} \leq \Alpha_k^{\min} + \theta \beta_k = (\Alpha' + \theta)\beta_k$, Lemma~\ref{lem.d_diff_bound}, Lemma~\ref{lem.gd_diff}, and $\beta_k \in (0,1]$ that
  \begin{align*}
    &\ \E_k[\phi(X_{k+1},\Tau_k) - \phi(X_k,\Tau_k)] = \E_k[\phi(X_k + \Alpha_k D_k,\Tau_k) - \phi(X_k,\Tau_k)] \\
    \leq&\ \E_k[-\Alpha_k \Delta l(X_k,\Tau_k,\nabla f(X_k),\Dtrue_k) \\
    &\hspace{0.6cm} + \Alpha_k \Tau_k \nabla f(X_k)^T(D_k - \Dtrue_k) + (1-\eta) \Alpha_k \beta_k \Delta l(X_k,\Tau_k,G_k,D_k)] \\
    \leq&\ - \Alpha' \beta_k \Delta l(X_k,\Tau',\nabla f(X_k),\Dtrue_k) \\
    &\ + (\Alpha' + \theta) \beta_k \Tau' \kappa_{\nabla f} \zeta^{-1} \sqrt{\rho_k} \\
    &\ + (1-\eta) (\Alpha' + \theta) \beta_k^2 (\Delta l(X_k, \Tau', \nabla f(X_k),\Dtrue_k) + \Tau' \zeta^{-1} (\rho_k + \kappa_{\nabla f} \sqrt{\rho_k})) \\
    \leq&\ - \Alpha' \beta_k \Delta l(X_k,\Tau',\nabla f(X_k),\Dtrue_k) \\
    &\ + (\Alpha' + \theta) \beta_k \Tau' \kappa_{\nabla f} \zeta^{-1} \sqrt{\iota} \beta_k \\
    &\ + (1-\eta) (\Alpha' + \theta) \beta_k^2 (\Delta l(X_k,\Tau',\nabla f(X_k),\Dtrue_k) + \Tau' \zeta^{-1} (\iota \beta_k^2 + \kappa_{\nabla f} \sqrt{\iota} \beta_k)) \\
    \leq&\ -\beta_k (\Alpha' - (1-\eta)(\Alpha' + \theta)\beta_k) \Delta l(X_k,\Tau',\nabla f(X_k),\Dtrue_k) + R' \beta_k^2,
  \end{align*}
  where $R' = (\Alpha' + \theta) \Tau' \zeta^{-1} (\kappa_{\nabla f} \sqrt{\iota} + (1-\eta) (\iota + \kappa_{\nabla f} \sqrt{\iota}))$.  Now, from Assumption~\ref{ass.good_merit_parameter} (which subsumes Assumption~\ref{ass.prob}), there exists $\phi_{\min} \in \R{}$ such that $\phi(X_k,\Tau') \geq \phi_{\min}$ for all $k \geq k_{\min}$.  One also finds that $\alpha_{\min}' \leq \Alpha' \leq \alpha_{\max}'$ due to the monotonicity of $\frac{2(1-\eta)\Xi'\tau}{\tau L + \Gamma}$ with respect to $\tau$.  Therefore, under part $(a)$ of the lemma, in which case one finds for $k \geq k_{\min}$ that $\psi \frac{\alpha_{\min}'}{2 (1-\eta) (\alpha_{\min}' + \theta)} \leq \beta \leq \psi \frac{\alpha_{\max}'}{2 (1-\eta) (\alpha_{\max}' + \theta)}$, it follows from above that
  \begin{align*}
    &\ \E_k[\phi(X_{k+1},\Tau_k) - \phi(X_k,\Tau_k)] \\
    \leq &\ - \frac{\psi (1 - \frac\psi2) (\alpha_{\min}')^2}{2 (1-\eta) (\alpha_{\min}' + \theta)} \Delta l(X_k,\Tau',\nabla f(X_k),\Dtrue_k) + \rho_{\max}' \(\psi \frac{\alpha_{\max}'}{2 (1-\eta) (\alpha_{\max}' + \theta)}\)^2
  \end{align*}
  so by taking total expectation conditioned on the event $\Ecal$ one finds
  \begin{align*}
    &\ \phi_{\min} - \E[\phi(X_{k_{\min}}, \Tau') | \Ecal ] \\
    \leq&\ \E[\phi(X_{k_{\min}+k}, \Tau') - \phi(X_{k_{\min}}, \Tau') | \Ecal] = \E\left[ \sum_{j = k_{\min}}^{k_{\min}+k-1} (\phi(X_{j+1}, \Tau') - \phi(X_j,\Tau')) \Bigg| \Ecal \right] \\
    \leq&\ - \tfrac{\psi (1 - \frac\psi2) (\alpha_{\min}')^2}{2 (1-\eta) (\alpha_{\min}' + \theta)} \E\left[ \sum_{j = k_{\min}}^{k_{\min}+k-1} \Delta l(X_j,\Tau',\nabla f(X_j),\Dtrue_j) \Bigg| \Ecal \right] + k \rho_{\max}' \(\psi \tfrac{\alpha_{\max}'}{2 (1-\eta) (\alpha_{\max}' + \theta)}\)^2.
  \end{align*}
  Rearranging terms, observing that $\E[\phi(X_{k_{\min}}, \Tau') | \Ecal ]$ is bounded above under Assumption~\ref{ass.good_merit_parameter}, and considering the limit superior as $k \to \infty$, the conclusion of part $(a)$ follows.  On the other hand, under the conditions of part $(b)$, it follows in a similar manner that, for any $k \in \N{}$, one finds
  \begin{align*}
    &\ \phi_{\min} - \E[\phi(X_{k_{\min}}, \Tau') | \Ecal ] \\
    \leq&\ \E[\phi(X_{k_{\min}+k}, \Tau') - \phi(X_{k_{\min}}, \Tau') | \Ecal ] = \E\left[\sum_{j = k_{\min}}^{k_{\min}+k-1} ( \phi(X_{j+1},\Tau') - \phi(X_j,\Tau')) \Bigg| \Ecal \right] \\
    \leq&\ \E\left[\sum_{j = k_{\min}}^{k_{\min}+k-1} (-\beta_j (\Alpha' - (1-\eta)(\Alpha' + \theta) \beta_j) \Delta l(X_j,\Tau',\nabla f(X_j),\Dtrue_j) + R' \beta_j^2 ) \Bigg| \Ecal \right].
  \end{align*}
  Taking limits as $k \to \infty$, the conclusion of part $(b)$ follows.
\eproof

We now present our main convergence theorem for Algorithm~\ref{alg.sqp}, which is essentially a translation of Lemma~\ref{lem.stoch_main} from results about model reductions to results about quantities connected to measures of stationarity for problem~\eqref{prob.opt}.

\btheorem\label{th.main}
  Suppose the conditions of Lemma~\ref{lem.stoch_main} hold.  Then,
  \benumerate
    \item[(a)] under the conditions of Lemma~\ref{lem.stoch_main}(a), there exists $C \in \R{}_{>0}$ such that
  \eenumerate
  \bequationNN
    \limsup_{k\to\infty}\ \E\left[\frac{1}{k} \sum_{j=k_{\min}}^{k_{\min}+k-1} (\thalf \Tau' \zeta \|\Dtrue_j\|_2^2 + \sigma(\|c(X_j)\|_2 - \|c(X_j)+J_j\Dtrue_j\|_2)) \Bigg| \Ecal \right] = C;
  \eequationNN
  \benumerate
    \item[(b)] under the conditions of Lemma~\ref{lem.stoch_main}(b), with $B_k := \sum_{j = k_{\min}}^{k_{\min} + k - 1} \beta_j$,
  \eenumerate
  \bequationNN
    \E\left[ \frac{1}{B_k} \sum_{j=k_{\min}}^{k_{\min} + k - 1} \beta_j (\thalf \Tau' \zeta \|\Dtrue_j\|_2^2 + \sigma (\|c(X_j)\|_2 - \|c(X_j) + J_j\Dtrue_j\|_2)) \Bigg| \Ecal \right] \xrightarrow{k \to \infty} 0,
  \eequationNN
  which further implies $\liminf_{k\to\infty}\ \E[\|\Dtrue_k\|_2^2 + (\|c(X_k)\|_2 - \|c(X_k) + J_k\Dtrue_k\|_2) | \Ecal] = 0$. 
\etheorem
\bproof
  The desired conclusions follow from Lemmas~\ref{lem.daniel} and \ref{lem.stoch_main}.
\eproof

One might be able to strengthen the conclusion in Theorem~\ref{th.main}(b), say to an almost-sure convergence guarantee; see, e.g., \cite{BertTsit00}.  However, we are satisfied with Theorem~\ref{th.main}(b), which is sufficient for revealing the favorable properties of Algorithm~\ref{alg.sqp} under Assumptions~\ref{ass.g},~\ref{ass.H}, and~\ref{ass.good_merit_parameter}.  Theorem~\ref{th.main}(a) shows under Assumptions~\ref{ass.g},~\ref{ass.H}, and~\ref{ass.good_merit_parameter} that if the latter condition in \eqref{eq.g} holds with $\rho_k = \rho$ for some $\rho \in \R{}_{>0}$ for all $k \in \N{}$ and $\{\beta_k\} = \{\beta\}$ is chosen as a (sufficiently small) constant sequence, then the limit superior of the expectation of the average of quantities connected to stationarity measures for problem~\eqref{prob.opt} is bounded above by a constant proportional to $\beta$. Intuitively, this shows that the iterates generated by the algorithm ultimately remain in a region in which these stationarity measures are small.  On the other hand, Theorem~\ref{th.main}(b) shows under Assumption~\ref{ass.good_merit_parameter} that if $\{\rho_k\}$ and $\{\beta_k\}$ vanish with $\rho_k = \Ocal(\beta_k^2)$, then a subsequence of iterates exist over which the expected values of these stationarity measures vanish.  As seen in Lemma~\ref{lem.normal_stationarity_2}, if there exists a subsequence of iterates, say indexed by $\Scal \subseteq \N{}$, that converges to a point satisfying certain regularity conditions, then $\{\|c_k\|_2 - \|c_k + \nabla c(x_k)^Tv_k\|_2\}_{k \in \Scal} \to 0$ means that the limit point is stationary with respect to the problem to minimize $\thalf \|c(x)\|_2^2$ subject to $x \in \R{n}_{\geq0}$.  Similarly, as seen in Lemma~\ref{lem.stationarity}, if there exists such a subsequence and the limit point is feasible with respect to problem~\eqref{prob.opt}, then $\{\dtrue_k\}_{k \in \Scal} \to 0$ means that the limit point is stationary with respect to \eqref{prob.opt}.  These situations are not guaranteed to occur, but this discussion shows that Theorem~\ref{th.main} is meaningful.

%************
% Subsection
%************
\subsection{Non-vanishing Merit Parameter}\label{sec.non-vanishing}

Our main convergence result in the previous section, namely, Theorem~\ref{th.main}, requires Assumption~\ref{ass.good_merit_parameter}, which in turn requires that the merit parameter sequence ultimately becomes a sufficiently small, positive constant sequence.  (Recall the discussion after Assumption~\ref{ass.good_merit_parameter}.)  To show that this corresponds to a realistic event for practical purposes, we next show conditions under which one finds that the merit parameter would not vanish.

We begin by showing a generally applicable result about the solution of \eqref{prob.v}.  It is related to that in Lemma~\ref{lem.normal_stationarity_2}, but is stronger due to an additional assumption.

\blemma\label{lem.wk_2_lb}
  Suppose the conditions of Lemma~\ref{lem.normal_stationarity_2} hold and  there exists $\kappa_w \in [0,1)$ such that for all generated $k \in \N{}$ in any run of the algorithm one has $\|c_k + \nabla c(x_k)^Tv_k\|_2 \leq \kappa_w \|c_k\|_2$.  Then, there exists $\kappa_v \in \R{}_{>0}$ such that, in any run of the algorithm such that iteration $k \in \N{}$ is reached, one finds
  \bequation\label{eq.v2}
    \|c_k\|_2 - \|c_k + \nabla c(x_k)^Tv_k\|_2 \geq \kappa_v \|v_k\|_2.
  \eequation
\elemma
\bproof
  Consider an arbitrary run of the algorithm in which the conditions of the lemma hold and iteration $k \in \N{}$ is reached.  If $c_k = 0$, then it follows by construction of \eqref{prob.v} that $v_k = 0$, in which case \eqref{eq.v2} follows trivially.  Hence, we may proceed under the assumption that $c_k \neq 0$, which by the conditions of the lemma, Assumption~\ref{ass.prob} (see \eqref{eq.bounds}), and the triangle inequality gives
  \bequationNN
    \kappa_{\nabla c} \|v_k\|_2 \geq \|\nabla c(x_k)^Tv_k\|_2 \geq \|c_k\|_2 - \|c_k + \nabla c(x_k)^Tv_k\|_2 \geq (1 - \kappa_w) \|c_k\|_2.
  \eequationNN
  Consequently, from \eqref{prob.alpha01}, \eqref{eq.intermediate}, and a similar derivation as in Lemma~\ref{lem.normal_stationarity_2}, one finds
  \begin{align*}
    &2\|c_k\|_2 (\|c_k\|_2 - \|c_k + \nabla c(x_k)^Tv_k\|_2) \geq \|c_k\|_2^2 - \|c_k + \nabla c(x_k)^Tv_k\|_2^2 \\
    \geq \ &\min\{\tfrac{\lambda^2}{\kappa_{\nabla c}^2}, 2\mu\} \|v_k\|_2^2 \geq \min\{\tfrac{\lambda^2}{\kappa_{\nabla c}^2}, 2\mu\} (\tfrac{1 - \kappa_w}{\kappa_{\nabla c}}) \|c_k\|_2 \|v_k\|_2,
  \end{align*}
  from which the desired conclusion in \eqref{eq.v2} follows.
\eproof

We now show that, under common conditions and when the norm of the stochastic gradient estimate is bounded uniformly, the denominator of the formula for $\tautrial_k$ in \eqref{eq.tautrial} is bounded proportionally to $\|v_k\|_2$.

\blemma\label{lem.tautrial_denominator}
  Suppose that Assumptions~\ref{ass.prob} and \ref{ass.H} hold, and that there exists $(\lambda,\mu,\kappa_g) \in \R{}_{>0} \times \R{}_{>0} \times \R{}_{>0}$ such that for all generated $k \in \N{}$ in any run of the algorithm one has $\nabla c(x_k)^T\nabla c(x_k) \succeq \lambda I$, $\mu_k \geq \mu$, and $\|g_k\|_2 \leq \kappa_g$.  Then, there exists $\kappa_{g,H} \in \R{}_{>0}$ such that, in any run such that iteration $k \in \N{}$ is reached, one finds
  \bequationNN
    g_k^Td_k + \thalf d_k^TH_kd_k \leq \kappa_{g,H} \|v_k\|_2.
  \eequationNN
\elemma
\bproof
  Consider an arbitrary run in which the conditions of the lemma hold and iteration $k \in \N{}$ is reached.  By Lemma~\ref{lem.well_defined}, $(u,w) = (0,0)$ is feasible for \eqref{prob.v}, so
  \begin{align*}
    &\ \max\left\{\thalf\|c_k + \nabla c(x_k)^T\nabla c(x_k)w_k\|_2^2, \thalf \mu_k \|u_k\|_2^2 \right\} \\
    \leq&\ \thalf\|c_k + \nabla c(x_k)^T\nabla c(x_k)w_k\|_2^2 + \thalf \mu_k \|u_k\|_2^2 \leq \thalf\|c_k\|_2^2.
  \end{align*}
  Since $\thalf\|c_k + \nabla c(x_k)^T\nabla c(x_k)w_k\|_2^2 \leq \thalf\|c_k\|_2^2$, it follows that
  \bequationNN
    \|\nabla c(x_k)^T\nabla c(x_k)w_k\|_2^2 \leq -2c_k^T\nabla c(x_k)^T\nabla c(x_k)w_k \leq 2\|c_k\|_2\|\nabla c(x_k)^T\nabla c(x_k)w_k\|_2,
   \eequationNN
   which along with Assumption~\ref{ass.prob} (see \eqref{eq.bounds}) shows that
  \bequationNN
    \|\nabla c(x_k)w_k\|_2 \leq \kappa_{\nabla c} \|w_k\|_2 \leq \tfrac{\kappa_{\nabla c}}{\lambda}\|\nabla c(x_k)^T\nabla c(x_k)w_k\|_2 \leq 2 \tfrac{\kappa_{\nabla c}}{\lambda} \|c_k\|_2 \leq 2 \tfrac{\kappa_{\nabla c}}{\lambda} \kappa_c.
  \eequationNN
  On the other hand, since $\thalf \mu_k \|u_k\|_2^2 \leq \thalf \|c_k\|_2^2$, it follows under Assumption~\ref{ass.prob} that $\|u_k\|_2 \leq \tfrac{1}{\sqrt{\mu_k}} \|c_k\|_2 \leq \tfrac{1}{\sqrt{\mu}} \kappa_c$.  Therefore, overall, it follows that
  \bequationNN
    \|v_k\|_2 = \sqrt{\|\nabla c(x_k)w_k\|_2^2 + \|u_k\|_2^2} \leq \(\sqrt{4 (\tfrac{\kappa_{\nabla c}}{\lambda})^2 + \tfrac{1}{\mu}}\) \kappa_c.
  \eequationNN
  Now, since $v_k = \nabla c(x_k)w_k + u_k$ is a feasible solution of \eqref{prob.d} while $d_k$ is the optimal solution of \eqref{prob.d}, it follows under the conditions of the lemma that
  \begin{align*}
    g_k^Td_k + \thalf d_k^TH_kd_k
    &\leq g_k^Tv_k + \thalf v_k^TH_kv_k \\
    &\leq \kappa_g \|v_k\|_2 + \thalf \kappa_H \|v_k\|_2^2 \leq \(\kappa_g + \thalf \kappa_H \(\sqrt{4 (\tfrac{\kappa_{\nabla c}}{\lambda})^2 + \tfrac{1}{\mu}}\) \kappa_c \) \|v_k\|_2,
  \end{align*}
  which leads to the desired conclusion.
\eproof

We now prove conditions under which the merit parameter does not vanish.

\btheorem\label{thm.tau_away_from_zero}
  Suppose that Assumptions~\ref{ass.prob} and \ref{ass.H} hold, and that there exists $(\lambda,\mu,\kappa_g,\kappa_w) \in \R{}_{>0} \times \R{}_{>0} \times \R{}_{>0} \times [0,1)$ such that for all generated $k \in \N{}$ in any run of the algorithm one has $\nabla c(x_k)^T\nabla c(x_k) \succeq \lambda I$, $\mu_k \geq \mu$, $\|g_k\|_2 \leq \kappa_g$, and $\|c_k + \nabla c(x_k)^Tv_k\|_2 \leq \kappa_w \|c_k\|_2$.  Then, in any run that does not terminate finitely, the latter event in Lemma~\ref{lem.well_defined}$(k)$ occurs $($i.e., $\{\tau_k\}$ does not vanish$)$ with $\tau_{\min} \geq \tfrac{(1-\sigma)\kappa_v}{\kappa_{g,H}} (1 - \epsilon_\tau)$.
\etheorem
\bproof
  Consider arbitrary $k \in \N{}$ in a run that does not terminate finitely and note that if $d_k = 0$ or $g_k^Td_k + \thalf d_k^TH_kd_k \leq 0$, then $\tautrial_k \gets \infty$, and otherwise $\tautrial_k$ is set by \eqref{eq.tautrial}.  Hence, under the conditions of the lemma and by Lemmas~\ref{lem.wk_2_lb}--\ref{lem.tautrial_denominator},
  \bequationNN
    \tautrial_k \geq \tfrac{(1-\sigma)(\|c_k\|_2 - \|c_k + \nabla c(x_k)^Td_k\|_2)}{g_k^Td_k + \thalf d_k^TH_kd_k} = \tfrac{(1-\sigma)(\|c_k\|_2 - \|c_k + \nabla c(x_k)^Tv_k\|_2)}{g_k^Td_k + \thalf d_k^TH_kd_k} \geq \tfrac{(1-\sigma)\kappa_v}{\kappa_{g,H}} =: \tau_*.
  \eequationNN
  Consequently, by the merit parameter update in \eqref{eq.merit_parameter}, $\tau_k < \tau_{k-1}$ only if $\tau_{k-1} > \tau_*$. This, along with Lemma~\ref{lem.well_defined}$(d)$--$(e)$, leads to the conclusion.
\eproof

Since $\nabla f$ is bounded in norm over the set $\Xcal$ in Assumption~\ref{ass.prob}, Theorem~\ref{thm.tau_away_from_zero} shows that, amongst the other stated conditions, if $\|g_k - \nabla f(x_k)\|_2$ is bounded uniformly over all $k \in \N{}$ in any, then the merit parameter sequence always remains bounded below by a positive number.  Under such conditions, the only potentially poor behavior of the merit parameter sequence is that, in a given run, it ultimately remains constant at a value that is too large.  We claim that, under certain assumptions about the distribution of the stochastic gradient estimates, this behavior can be shown to occur with probability zero.  (We do not prove such a result here, but refer the interested reader to Proposition 3.16 in \cite{BeraCurtRobiZhou21} to see such a result for the equality-constraints-only setting, in which case the behavior of the merit parameter is similar.)  On the other hand, if $\|g_k - \nabla f(x_k)\|_2$ is not bounded uniformly in this manner, then it is possible for the merit parameter sequence to vanish unnecessarily.  This issue is one that should be noted by a user of the algorithm.  In particular, if in a run of the algorithm one chooses $\mu_k \geq \mu$ for some $\mu \in \R{}_{>0}$ for all $k \in \N{}$ and one finds for some $(\lambda,\kappa_w) \in \R{}_{>0} \times \R{}_{>0}$ that generated $k \in \N{}$ yield $\nabla c(x_k)^T\nabla c(x_k) \succeq \lambda I$ and $\|c_k + \nabla c(x_k)^Tv_k\|_2 \leq \kappa_w \|c_k\|_2$, yet $\tau_k$ has become exceedingly small, then Theorem~\ref{thm.tau_away_from_zero} shows that this must be due to the stochastic gradient estimates tending to become significantly large in norm, in which case the performance of the algorithm may improve with more accurate stochastic gradient estimates.

%************
% Subsection
%************
\subsection{Deterministic Algorithm}\label{sec.deterministic}

We conclude this section with a statement of a convergence result that we claim to hold for Algorithm~\ref{alg.sqp} if it were to be run with $g_k = \nabla f(x_k)$ for all $k \in \N{}$.  Due to space considerations, we do not provide a proof of the result, although we offer the proposition for reference for the reader and claim that it holds from results proved in this paper for the stochastic setting as well as other similar results for SQP methods for deterministic continuous nonlinear optimization.

\begin{proposition}
  Suppose Assumptions~\ref{ass.prob} and~\ref{ass.H} hold and Algorithm~\ref{alg.sqp} is run with $g_k = \nabla f(x_k)$ for all $k \in \N{}$.  If for all large $k \in \N{}$ there exists $\kappa_w \in [0,1)$ such that $\|c_k + \nabla c(x_k)^Tv_k\|_2 \leq \kappa_w \|c_k\|_2$, then $\{x_k\} \subset \R{n}_{\geq0}$, $\{\tau_k\}$ is bounded away from zero, and, with $y_k \in \R{m}$ and $z_k \in \R{n}_{\geq0}$ defined as the optimal multipliers corresponding to the solution of subproblem~\eqref{prob.d} for all $k \in \N{}$, it follows that
  \bequationNN
    \left\{\left\|\bbmatrix \nabla f(x_k) + \nabla c(x_k) y_k - z_k \\ c_k \\ x_k^Tz_k \ebmatrix \right\| \right\} \to 0.
  \eequationNN
  Otherwise, $\{x_k\} \subset \R{n}_{\geq0}$, $\{\min\{\nabla c(x_k)c_k,0\}\} \to 0$, and $\{|x_k^T\nabla c(x_k)c_k|\} \to 0$, and if the sequence $\{\tau_k\}$ is bounded away from zero, then
  \bequationNN
    \left\{\left\|\bbmatrix \nabla f(x_k) + \nabla c(x_k) y_k - z_k \\ x_k^Tz_k \ebmatrix \right\| \right\} \to 0.
  \eequationNN
\end{proposition}

%*********
% Section
%*********
\section{Numerical Results}\label{sec.numerical}

In this section, we provide results demonstrating the performance of a MATLAB implementation of Algorithm~\ref{alg.sqp} when solving a subset of problems from CUTEst \cite{GoulOrbaToin15}, where Gurobi is used to solve the arising subproblems \cite{Gurobi23}.  The purpose of these experiments is to compare this performance against that of the Julia implementation provided by the authors of \cite[Algorithm 1]{NaAnitKola21}.  From all inequality-constrained problems in CUTEst, we selected those such that (i) $m \leq n \leq 1000$, (ii) $f(x_k) \geq -10^{20}$ for all $k \in \N{}$ in all runs of our algorithm, and (iii) Gurobi did not report any errors.  This resulted in a set of 323 test problems.

For each test problem, both codes used the same initial iterate and generated stochastic gradient estimates in the same manner.  Specifically, for all $k\in\N{}$ in each run, the codes set $g_k = \Ncal(\nabla f(x_k), \epsilon_g(I + ee^T))$, where $e$ is the all-ones vector and $\epsilon_g \in \{10^{-8},10^{-4},10^{-2},10^{-1}\}$ was fixed for each run (see below).  If a problem had only inequality constraints, i.e., $m=0$, then our code explicitly computed $\alpha_k^{\varphi}$ (as defined in \eqref{eq.alpha_max}) and set $\alpha_k \gets \alpha_k^{\max}$ for all $k \in \N{}$. Otherwise, the code set $\alpha_k \gets \min \{1,(1.1)^{t_k}\alpha_k^{\min},\alpha_k^{\min} + \theta\beta_k\}$, where $t_k \gets \max\{t \in \N{}: \varphi_k((1.1)^t\alpha_k^{\min}) \leq 0\}$.  This guarantees that $\alpha_k\in[\alpha_k^{\min},\alpha_k^{\max}]$ for all $k\in\N{}$.  The other user-defined parameters of Algorithm~\ref{alg.sqp} were selected as $\sigma = \tau_0 = 0.1$, $\eta = 0.5$, $\xi_0 = 1$, $\epsilon_{\tau} = \epsilon_{\xi} = 10^{-2}$, $\theta = 10^4$, $\mu_k = \max\{10^{-8},10^{-4}\|c_k\|_2^2\}$, $\beta_k = 1$, and $H_k = I$ for all $k\in\N{}$.  The Lipschitz constants $L$ and $\Gamma$ were estimated every 100 iterations by differences of stochastic gradients at ten samples around the current iterate.  Meanwhile, we ran the Julia code for \cite[Algorithm 1]{NaAnitKola21} with the \texttt{AdapGD} option and its default parameter settings as described in \cite[Section 4]{NaAnitKola21}.  Each code terminated as soon as $10^4$ stochastic gradient samples were evaluated or a 12-hour CPU time limit was reached.

Let $\texttt{FeasErr}(x)$ be the $\infty$-norm constraint violation at $x$ and let $\texttt{KKTErr}(x,y,z)$ be the $\infty$-norm violation of the KKT conditions (recall \eqref{eq.KKT}) at a primal-dual iterate $(x,y,z)$.  Each run of Algorithm~\ref{alg.sqp} generates $\{x_k\}\subset\R{n}$.  For each $k \in \N{}$, let $y_k^{\rm true}\in\R{m}$ and $z_k^{\rm true}\in\R{n}$ denote the optimal Lagrange multipliers corresponding to the equality and inequality constraints when \eqref{prob.d} is solved with $g_k = \nabla f(x_k)$.  For each run of Algorithm~\ref{alg.sqp}, we determined the best iterate as $x_{k_\texttt{best}}$ where
\begin{equation*}
k_{\texttt{best}} = \begin{cases} \displaystyle\arg\min_{k\in\N{}} \ \texttt{FeasErr}(x_k) \quad\quad\quad\quad\quad\ \ \text{if }\texttt{FeasErr}(x_k) > 10^{-4} \text{ for all } k\in\N{}, \\
\displaystyle\arg\min_{k\in\N{}}\
 \{\texttt{KKTErr}(x_k,y_k^{{\rm true}},z_k^{\rm true}): \texttt{FeasErr}(x_k) \leq 10^{-4}\} \quad\text{ otherwise.} \end{cases}
\end{equation*}
We determined the best iterate in a run of \cite[Algorithm 1]{NaAnitKola21} using the same formula with the sequence of iterates and Lagrange multiplier estimates that are computed as part of the algorithm. Our results for four noise levels, provided in Figure~\ref{fig.perf_comparison} below, are presented in terms of $\texttt{FeasErr}(x_{k_\texttt{best}})$ as the feasibility error and $\texttt{KKTErr}(x_{k_\texttt{best}},y_{{k_\texttt{best}}}^{\rm true},z_{{k_\texttt{best}}}^{\rm true})$ as the KKT error for each run of both algorithms.

Since the Julia code for \cite[Algorithm 1]{NaAnitKola21} is only set up to solve CUTEst problems without simple bound constraints, the results in Figure~\ref{fig.perf_comparison} are presented in two parts.  For the 57 problems for which both algorithms were set up to run, the first two box plots show the best feasibility and KKT errors achieved by both codes, where each problem is run 5 times each (since the behaviors of the algorithms are stochastic).  In the third box plot, we report the best feasibility and KKT errors obtained by our Matlab code on the remaining $266$ $(=323-57)$ problems, again with five runs for each problem.  Overall, one finds that the performance of our algorithm is comparatively good in this experimental set-up.  The best feasibility and KKT errors are relatively low for our algorithm, although the errors increase with the noise level, as may be expected.  Experiments with diminishing step sizes also showed favorable performance for our algorithm; these results are omitted due to page limit restrictions.

\begin{figure}[ht]
  \centering
 \includegraphics[width=0.28\textwidth,clip=true,trim=20 5 90 50]{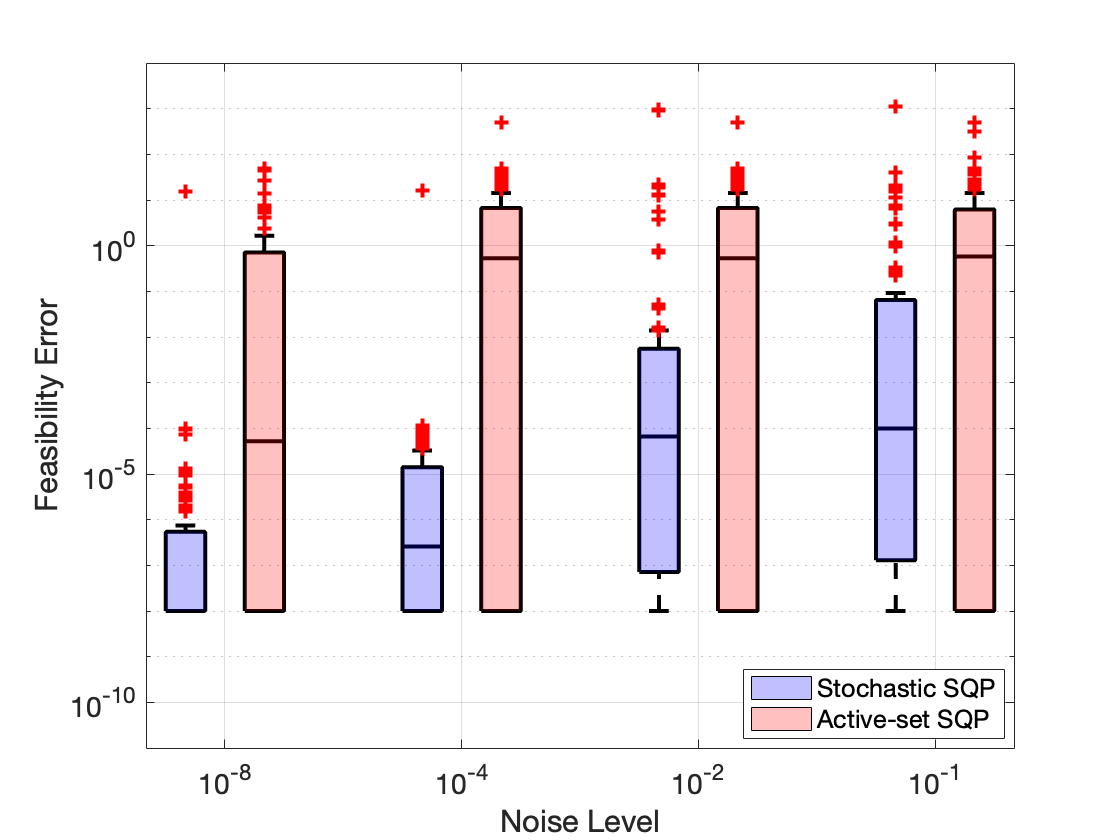}\qquad 
 \includegraphics[width=0.28\textwidth,clip=true,trim=20 5 90 50]{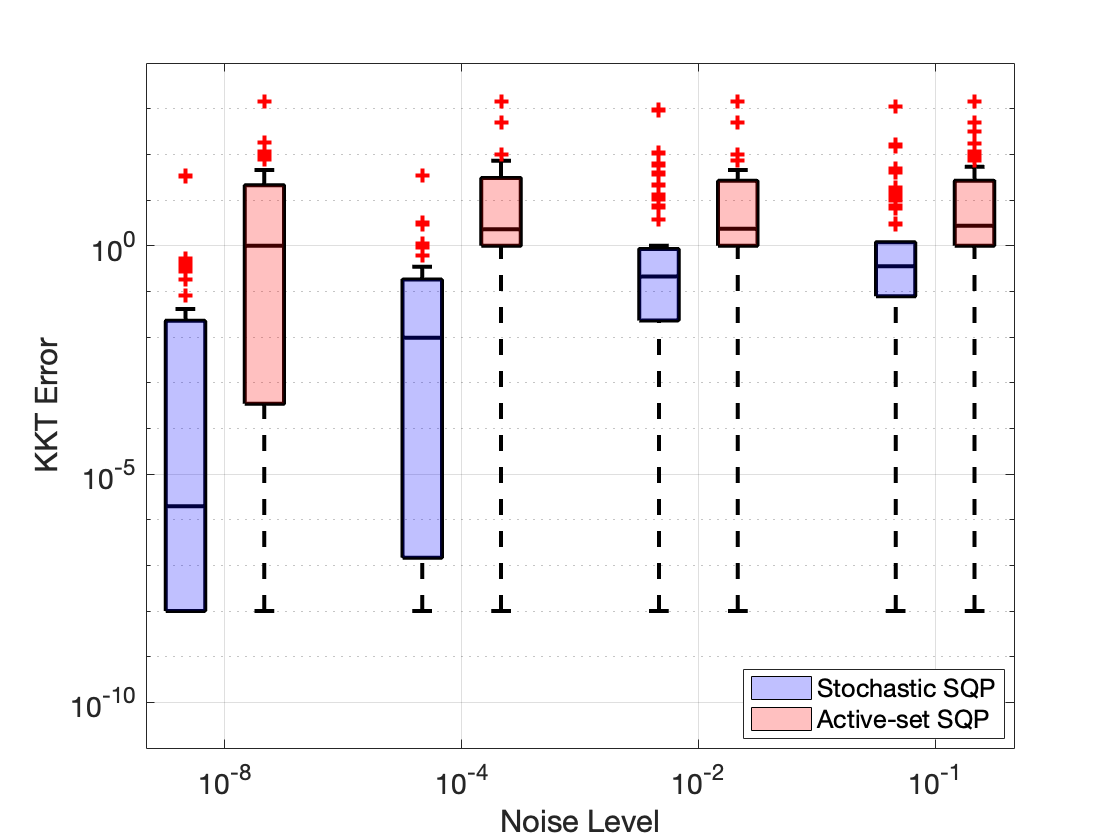}\qquad
 \includegraphics[width=0.28\textwidth,clip=true,trim=20 5 90 50]{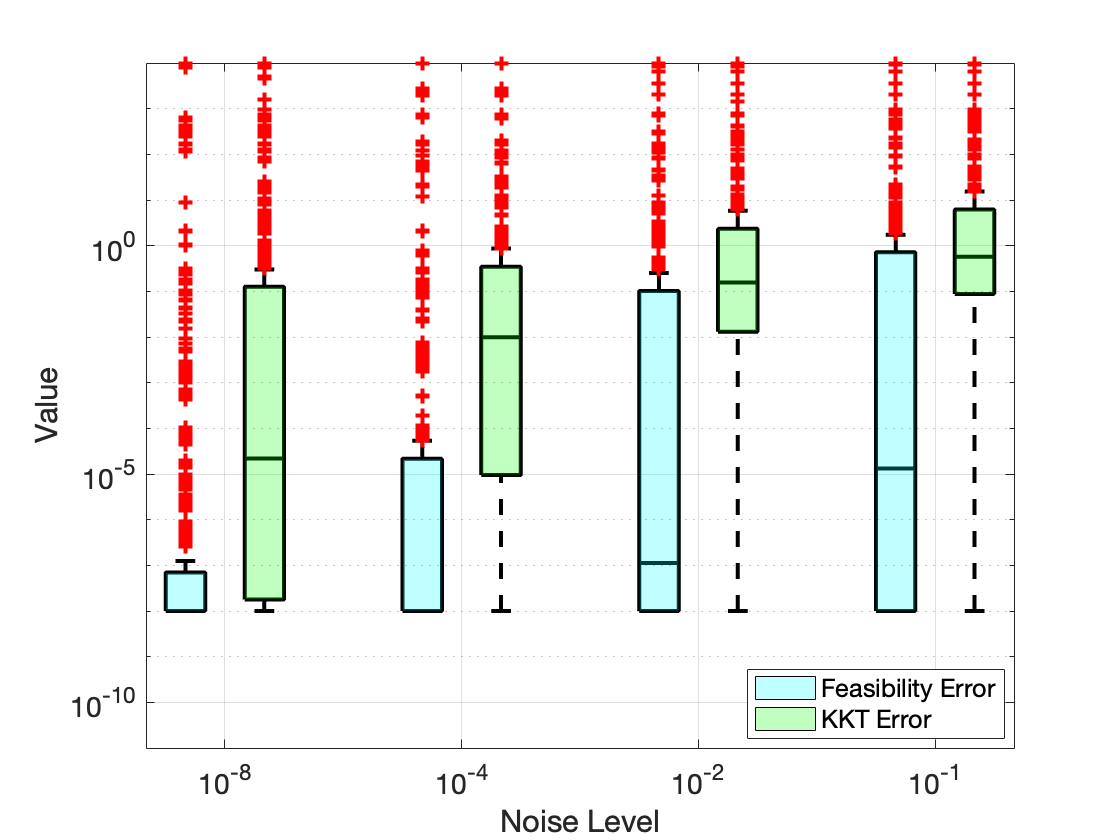}
  \caption{Box plots comparing the best feasibility errors (left) and KKT errors (middle) of a Matlab implementation of Algorithm~\ref{alg.sqp} (``Stochastic SQP") and the Julia implementation provided by the authors of \cite[Algorithm 1]{NaAnitKola21} (``Active-st SQP'') when solving 57 CUTEst problems without simple bound constraints.  Box plots of the best feasibility and KKT errors (combined, right) of the implementation of Algorithm~\ref{alg.sqp} when solving the other 266 CUTEst problems from the test set.}
  \label{fig.perf_comparison}
\end{figure}

% \begin{figure}[ht]
%   \centering
%   \includegraphics[width=0.35\textwidth,clip=true,trim=20 5 90 50]{our_method_infeasibility.png}\qquad 
%  \includegraphics[width=0.35\textwidth,clip=true,trim=20 5 90 50]{our_method_stationarity.png}
%   \caption{Performance of our method on CUTEst problems for feasibility (left) and KKT (right) errors.}
%   \label{fig.perf_our_method}
% \end{figure}

%*********
% Section
%*********
\section{Conclusion}\label{sec.conclusion}

We have proposed, analyzed, and tested an algorithm for solving continuous optimization problems.  The algorithm requires that constraint function and derivative values can be computed in each iteration, but does not require exact objective function and derivative values; rather, the algorithm merely requires that a stochastic objective gradient estimate is computed to satisfy relatively loose assumptions in each iteration.  The theoretical convergence guarantees of the algorithm require knowledge of Lipschitz constants for the objective gradient and constraint Jacobian, although in practice these constants can be estimated.  Our numerical experiments show that our proposed algorithm can outperform an alternative algorithm that relies on the ability to compute more accurate gradient estimates.  We have provided comments throughout the paper on how the assumptions that are required for our theoretical convergence guarantees might be loosened further.

%*********
% Section
%*********
\section*{Acknowledgements}

The authors are grateful to Sen Na for providing consultation about the Julia implementation provided by the authors of \cite[Algorithm 1]{NaAnitKola21}.  This material is based upon work supported by the U.S.~NSF under award CCF-2139735 and by the Office of Naval Research under award N00014-21-1-2532.

%**************
% Bibliography
%**************
\bibliographystyle{plain}
\bibliography{reference}

\begin{thebibliography}{10}

\bibitem{BeraCurtRobiZhou21}
A.~S. Berahas, F.~E. Curtis, D.~Robinson, and B.~Zhou.
\newblock Sequential quadratic optimization for nonlinear equality constrained
  stochastic optimization.
\newblock {\em SIAM Journal on Optimization}, 31(2):1352--1379, 2021.

\bibitem{BeraBollZhou22}
Albert~S Berahas, Raghu Bollapragada, and Baoyu Zhou.
\newblock An adaptive sampling sequential quadratic programming method for
  equality constrained stochastic optimization.
\newblock {\em arXiv preprint arXiv:2206.00712}, 2022.

\bibitem{BeraCurtOneiRobi21}
Albert~S Berahas, Frank~E Curtis, Michael~J O'Neill, and Daniel~P Robinson.
\newblock A stochastic sequential quadratic optimization algorithm for
  nonlinear equality constrained optimization with rank-deficient jacobians.
\newblock {\em arXiv preprint arXiv:2106.13015}, 2021.

\bibitem{BeraShiYiZhou22}
Albert~S Berahas, Jiahao Shi, Zihong Yi, and Baoyu Zhou.
\newblock Accelerating stochastic sequential quadratic programming for equality
  constrained optimization using predictive variance reduction.
\newblock {\em arXiv preprint arXiv:2204.04161}, 2022.

\bibitem{BeraXieZhou23}
Albert~S Berahas, Miaolan Xie, and Baoyu Zhou.
\newblock A sequential quadratic programming method with high probability
  complexity bounds for nonlinear equality constrained stochastic optimization.
\newblock {\em arXiv preprint arXiv:2301.00477}, 2023.

\bibitem{Bert09}
Dimitri Bertsekas.
\newblock {\em Convex Optimization Theory}, volume~1.
\newblock Athena Scientific, 2009.

\bibitem{Bert98}
Dimitri~P. Bertsekas.
\newblock {\em Network optimization: continuous and discrete models}, volume~8.
\newblock Athena Scientific, 1998.

\bibitem{BertTsit00}
Dimitri~P. Bertsekas and John~N. Tsitsiklis.
\newblock Gradient convergence in gradient methods with errors.
\newblock {\em SIAM Journal on Optimization}, 10(3):627--642, 2000.

\bibitem{ByrdGilbNoce00}
Richard~H Byrd, Jean~Charles Gilbert, and Jorge Nocedal.
\newblock A trust region method based on interior point techniques for
  nonlinear programming.
\newblock {\em Math. Prog.}, 89(1):149--185, 2000.

\bibitem{ByrdHribNoce99}
Richard~H Byrd, Mary~E Hribar, and Jorge Nocedal.
\newblock An interior point algorithm for large-scale nonlinear programming.
\newblock {\em SIAM Journal on Optimization}, 9(4):877--900, 1999.

\bibitem{Conn73}
Andrew~R Conn.
\newblock Constrained optimization using a nondifferentiable penalty function.
\newblock {\em SIAM Journal on Numerical Analysis}, 10(4):760--784, 1973.

\bibitem{CurtRobiZhou21}
F.~E. Curtis, D.~P. Robinson, and B.~Zhou.
\newblock Inexact sequential quadratic optimization for minimizing a stochastic
  objective function subject to deterministic nonlinear equality constraints.
\newblock {\em arXiv preprint arXiv:2107.03512}, 2021.

\bibitem{CurtOneiRobi21}
Frank~E Curtis, Michael~J O'Neill, and Daniel~P Robinson.
\newblock Worst-case complexity of an sqp method for nonlinear equality
  constrained stochastic optimization.
\newblock {\em arXiv preprint arXiv:2112.14799}, 2021.

\bibitem{DipiGrip89}
G~Di~Pillo and L~Grippo.
\newblock Exact penalty functions in constrained optimization.
\newblock {\em SIAM Journal on control and optimization}, 27(6):1333--1360,
  1989.

\bibitem{DipiGrip85}
Gianni Di~Pillo and Luigi Grippo.
\newblock A continuously differentiable exact penalty function for nonlinear
  programming problems with inequality constraints.
\newblock {\em SIAM Journal on Control and Optimization}, 23(1):72--84, 1985.

\bibitem{Diki67}
II~Dikin.
\newblock Iterative solution of problems of linear and quadratic programming.
\newblock In {\em Doklady Akademii Nauk}, volume 174, pages 747--748. Russian
  Academy of Sciences, 1967.

\bibitem{FangNaMahoKola22}
Yuchen Fang, Sen Na, Michael~W. Mahoney, and Mladen Kolar.
\newblock Fully stochastic trust-region sequential quadratic programming for
  equality-constrained optimization problems.
\newblock {\em arXiv preprint 2211.15943}, 2022.

\bibitem{Flet73}
Roger Fletcher.
\newblock An exact penalty function for nonlinear programming with
  inequalities.
\newblock {\em Mathematical Programming}, 5(1):129--150, 1973.

\bibitem{Flet13}
Roger Fletcher.
\newblock {\em Practical methods of optimization}.
\newblock John Wiley \& Sons, 2013.

\bibitem{GillMurrSaun05}
Philip~E Gill, Walter Murray, and Michael~A Saunders.
\newblock Snopt: An {SQP} algorithm for large-scale constrained optimization.
\newblock {\em SIAM review}, 47(1):99--131, 2005.

\bibitem{GladPola79}
Torkel Glad and Elijah Polak.
\newblock A multiplier method with automatic limitation of penalty growth.
\newblock {\em Mathematical Programming}, 17(1):140--155, 1979.

\bibitem{GoulOrbaToin15}
Nicholas~IM Gould, Dominique Orban, and Philippe~L Toint.
\newblock {CUTEst: a constrained and unconstrained testing environment with
  safe threads for mathematical optimization}.
\newblock {\em Computational optimization and applications}, 60(3):545--557,
  2015.

\bibitem{Gurobi23}
{Gurobi Optimization, LLC}.
\newblock {Gurobi Optimizer Reference Manual}, 2023.

\bibitem{Han76}
Shih-Ping Han.
\newblock Superlinearly convergent variable metric algorithms for general
  nonlinear programming problems.
\newblock {\em Mathematical Programming}, 11(1):263--282, 1976.

\bibitem{Hath85}
Richard~J Hathaway.
\newblock A constrained formulation of maximum-likelihood estimation for normal
  mixture distributions.
\newblock {\em The Annals of Statistics}, 13(2):795--800, 1985.

\bibitem{Ibar88}
Toshihide Ibaraki and Naoki Katoh.
\newblock {\em Resource allocation problems: algorithmic approaches}.
\newblock MIT press, 1988.

\bibitem{Kour16}
Drew~P Kouri and Thomas~M Surowiec.
\newblock Risk-averse {PDE}-constrained optimization using the conditional
  value-at-risk.
\newblock {\em SIAM Journal on Optimization}, 26(1):365--396, 2016.

\bibitem{Lan20}
Guanghui Lan.
\newblock {\em First-order and stochastic optimization methods for machine
  learning}.
\newblock Springer, 2020.

\bibitem{LasdWareRice67}
L~Lasdon, A~Waren, and R~Rice.
\newblock An interior penalty method for inequality constrained optimal control
  problems.
\newblock {\em IEEE Transactions on Automatic Control}, 12(4):388--395, 1967.

\bibitem{Mcgi65}
Robert McGill.
\newblock Optimum control, inequality state constraints, and the generalized
  newton-raphson algorithm.
\newblock {\em Journal of the Society for Industrial and Applied Mathematics,
  Series A: Control}, 3(2):291--298, 1965.

\bibitem{NaAnitKola21}
Sen Na, Mihai Anitescu, and Mladen Kolar.
\newblock Inequality constrained stochastic nonlinear optimization via
  active-set sequential quadratic programming.
\newblock {\em arXiv preprint arXiv:2109.11502}, 2021.

\bibitem{NaAnitKola22}
Sen Na, Mihai Anitescu, and Mladen Kolar.
\newblock An adaptive stochastic sequential quadratic programming with
  differentiable exact augmented lagrangians.
\newblock {\em Mathematical Programming}, pages 1--71, 2022.

\bibitem{NaMaho22}
Sen Na and Michael~W Mahoney.
\newblock Asymptotic convergence rate and statistical inference for stochastic
  sequential quadratic programming.
\newblock {\em arXiv preprint arXiv:2205.13687}, 2022.

\bibitem{NoceWrig06}
Jorge Nocedal and Stephen Wright.
\newblock {\em Numerical Optimization}.
\newblock Springer, 2006.

\bibitem{OztoByrdNoce21}
Figen Oztoprak, Richard Byrd, and Jorge Nocedal.
\newblock Constrained optimization in the presence of noise.
\newblock {\em arXiv preprint arXiv:2110.04355}, 2021.

\bibitem{PateZhan21}
Vivak Patel and Shushu Zhang.
\newblock Stochastic gradient descent on nonconvex functions with general noise
  models.
\newblock {\em arXiv preprint 2104.00423}, 2021.

\bibitem{Pero84}
Andre~F Perold.
\newblock Large-scale portfolio optimization.
\newblock {\em Management science}, 30(10):1143--1160, 1984.

\bibitem{QiuKung23}
Songqiang Qiu and Vyacheslav Kungurtsev.
\newblock A sequential quadratic programming method for optimization with
  stochastic objective functions, deterministic inequality constraints and
  robust subproblems.
\newblock {\em arXiv preprint arXiv:2302.07947}, 2023.

\bibitem{Robi07}
Daniel~P. Robinson.
\newblock {\em Primal-dual methods for nonlinear optimization (thesis)}.
\newblock University of California, San Diego, 2007.

\bibitem{ShapDentRusz21}
Alexander Shapiro, Darinka Dentcheva, and Andrzej Ruszczynski.
\newblock {\em Lectures on stochastic programming: modeling and theory}.
\newblock SIAM, 2021.

\bibitem{ShiWangWang22}
Qiankun Shi, Xiao Wang, and Hao Wang.
\newblock A momentum-based linearized augmented lagrangian method for nonconvex
  constrained stochastic optimization.
\newblock 2022.

\bibitem{WaecBieg06}
A.~W{\"a}chter and L.~T. Biegler.
\newblock On the implementation of an interior-point filter line-search
  algorithm for large-scale nonlinear programming.
\newblock {\em Mathematical Programming}, 106(1):25--57, 2006.

\bibitem{Wrig97}
Stephen~J Wright.
\newblock {\em Primal-dual interior-point methods}.
\newblock SIAM, 1997.

\bibitem{Yama98}
Hiroshi Yamashita.
\newblock A globally convergent primal-dual interior point method for
  constrained optimization.
\newblock {\em Optimization Methods and Software}, 10(2):443--469, 1998.

\bibitem{Zang67}
Willard~I Zangwill.
\newblock Non-linear programming via penalty functions.
\newblock {\em Management science}, 13(5):344--358, 1967.

\bibitem{ZavaAnit14}
Victor~M. Zavala and Mihai Anitescu.
\newblock Scalable nonlinear programming via exact differentiable penalty
  functions and trust-region newton methods.
\newblock {\em SIAM Journal on Optimization}, 24(1):528--558, 2014.

\end{thebibliography}

\end{document}